\DocumentMetadata{lang=en-US}
\documentclass[
	a4paper,
	11pt,
	abstract=on, 
	numbers=auto, 
	toc=bib, 
	toc=numberline 
]{scrartcl}

\usepackage[vmargin={0.8in},
	hmargin={0.8in},
	footskip=0.4in, 
	nomarginpar, %
]{geometry} 
\NewDocumentCommand{\blfootnote}{O{}m}{%
	\begingroup
	\protect\renewcommand{\thefootnote}{#1}\footnote{#2}%
	\addtocounter{footnote}{-1}%
	\endgroup
}
\usepackage{xspace} 
\usepackage{enumitem} 
\setlist{noitemsep,topsep=0pt} 
\usepackage{array} 
\usepackage[export]{adjustbox} 
\usepackage[font=small]{caption} 
\usepackage{subcaption} 
\usepackage{float} 
\usepackage[svgnames]{xcolor} 
\usepackage{pifont} 
\usepackage{iftex} 
\usepackage{babel} 
\ifpdftex
	\usepackage[utf8]{inputenc}
\fi

\usepackage{amsmath}
\usepackage{mathtools}
\usepackage{amsthm}

\ifpdftex
	\usepackage[T1]{fontenc} 

	\usepackage[proportional,lining,subscriptcorrection,varg,varvw,varbb]{newtx}
	\usepackage[defaultsans,scale=0.95]{opensans}
\fi
\ifluatex
	\usepackage{fontspec} 
	\defaultfontfeatures{Ligatures=TeX,Scale=MatchLowercase,Numbers={Proportional,Lining}} 
	\usepackage{unicode-math} 

\fi
\usepackage{microtype} 

\setkomafont{disposition}{\bfseries\sffamily}
\setkomafont{title}{\bfseries\sffamily}
\setkomafont{author}{\normalsize}
\setkomafont{date}{\normalsize}
\makeatletter
\xpatchcmd{\@maketitle}{\huge}{\LARGE}{}{}
\makeatother
\addtokomafont{section}{\Large}
\addtokomafont{subsection}{\large}

\usepackage{tabularray} 
\usepackage{extarrows}
\usepackage{tikz}
\usetikzlibrary{arrows.meta,matrix,decorations.pathreplacing,shapes.geometric,backgrounds,calc,math}
\usepackage{tikz-cd}
\usepackage{makecell} 
\usepackage{algorithm2e}

\usepackage{hyperref} 
\usepackage{nameref} 
\usepackage[nameinlink]{cleveref} 
\ifdefined\pdfinfoomitdate
\fi
\ifdefined\pdfsuppressptexinfo
\fi
\ifdefined\pdfvariable
	\pdfvariable suppressoptionalinfo 1023 \relax
\fi
\makeatletter
\NewDocumentCommand{\keywords}{m}{\def\@keywords{#1}}
\AtBeginDocument{%
	\hypersetup{
		pdftitle = {\@title},
		pdfauthor = {\@author},
		pdfsubject = {\@subject},
		pdfkeywords = {\@keywords},
		pdfcreator = {Not T3X},
		pdfproducer = {Definitely not T3X}
	}
}
\makeatother
\newcommand*{\fullref}[1]{\hyperref[{#1}]{\nameref*{#1} (\ref*{#1})}}

\newcommand*{\locallabel}[1]{\csname label\endcsname{locallabel: sec: \thesection subsec: \thesubsection thm: \thetheorem label: #1}}
\newcommand*{\localref}[1]{\hyperref[locallabel: sec: \thesection subsec: \thesubsection thm: \thetheorem label: #1]{(\ref*{locallabel: sec: \thesection subsec: \thesubsection thm: \thetheorem label: #1})}}

\newtheorem*{theorem*}{Theorem}

\newtheorem*{lemma*}{Lemma}

\newtheorem*{proposition*}{Proposition}

\newtheorem*{corollary*}{Corollary}
\theoremstyle{definition}

\newtheorem*{definition*}{Definition}

\newtheorem*{example*}{Example}

\newtheorem*{remark*}{Remark}

\newcounter{diagram}
\makeatletter
\ExplSyntaxOn
\NewDocumentEnvironment{diagram*}{o}
{
\begin{equation*}
\IfNoValueTF{#1}{\begin{tikzcd}}{\begin{tikzcd}[#1]}
		}{
	\end{tikzcd}
	\end{equation*}
}
\NewDocumentEnvironment{diagram}{o}
{
\setcounter{diagram}{\value{equation}}
\stepcounter{equation}
\refstepcounter{diagram}
\begin{equation*}
\IfNoValueTF{#1}{\begin{tikzcd}}{\begin{tikzcd}[#1]}
		}{
	\end{tikzcd}
	\eqno\hbox{\normalfont\normalcolor(\thediagram)}
	\end{equation*}\@ignoretrue
}
\ExplSyntaxOff
\makeatother
\crefname{diagram}{diagram}{diagrams}
\Crefname{diagram}{Diagram}{Diagrams}
\creflabelformat{diagram}{#2(#1)#3}

\crefname{maintheorem}{theorem}{theorems}
\Crefname{maintheorem}{Theorem}{Theorems}

\crefname{mainproposition}{proposition}{propositions}
\Crefname{mainproposition}{Proposition}{Propositions}

\crefname{maincorollary}{corollary}{corollaries}
\Crefname{maincorollary}{Corollary}{Corollaries}

\usepackage[
	backend=biber,
	style=numeric-comp,
	sorting=none,
	isbn=true,
	url=false
]{biblatex}
\DeclareFieldFormat[misc,online]{title}{\mkbibquote{#1}}
\AtNextBibliography{\small}
\AtBeginDocument{
	\ifpdftex

	\fi
	
	\NewDocumentCommand{\ZZ}{}{\ifmmode\mathbb{Z}\fi}
	\NewDocumentCommand{\FF}{}{\ifmmode\mathbb{F}\fi}
	\NewDocumentCommand{\checkC}{}{\ifluatex Č\else\ifmmode\check{C}\else\v{C}\fi\fi} 
	\NewDocumentCommand{\Cech}{}{%
		\ifmmode%
			\operatorname{\checkC}%
		\else%
			\checkC ech\xspace\fi%
	} 
	\NewDocumentCommand{\DelCech}{}{\ifmmode\operatorname{ChrDel\checkC}\else Delaunay--\checkC ech\xspace\fi} 
}

\newcommand*{\Pipeline}{\texorpdfstring{M\textsuperscript{2}S\textsuperscript{2}}{M2S2}}
\newcommand*{\PipelineText}{M2S2}

\newcommand{\drawVectorRepresentation}[3]{%
	\def\squareSize{#1}
	\def\numSquares{#2}
	\def\rectHeight{\numSquares * \squareSize}

		\node[draw,rectangle,minimum width=\squareSize,minimum height=\rectHeight, anchor=north west] (#3) at (0,0) {};

		\foreach \i in {1,...,\numSquares} {
				\draw (0, -\i * \squareSize) -- (\squareSize, -\i * \squareSize);
			}
}

\title{Topology of Multi-species Localization}
\author{Abhinav Natarajan\textsuperscript{1}\and Thomas Chaplin\textsuperscript{1}\and Joshua A. Bull\textsuperscript{2}\and Eoghan J. Mulholland-Illingworth\textsuperscript{3}\and Simon J. Leedham\textsuperscript{3,5}\and Helen M. Byrne\textsuperscript{2,3,4}\and Maria-Jose Jimenez\textsuperscript{6,10}\and Heather A. Harrington\textsuperscript{7,8,9,10}}

\begin{document}
\maketitle
\footnotetext[1]{Mathematical Institute, University of Oxford, Radcliffe Observatory Quarter, Andrew Wiles Building, Woodstock Rd, Oxford OX2 6GG, UK.}
\footnotetext[2]{Wolfson Centre for Mathematical Biology, Mathematical Institute, University of Oxford, Oxford OX2 6GG, UK.}
\footnotetext[3]{Centre for Human Genetics, Roosevelt Drive, University of Oxford, Oxford, UK.}
\footnotetext[4]{Ludwig Institute for Cancer Research, Nuffield Department of Medicine, University of Oxford, Oxford OX3 7DQ, UK.}
\footnotetext[5]{Translational Gastroenterology Unit, John Radcliffe Hospital, University of Oxford, and Oxford NIHR Biomedical Research Centre, Oxford, UK.}
\footnotetext[6]{Universidad de Sevilla. Escuela Técnica Superior de Ingeniería Informática. Av. Reina Mercedes s/n, 41012 Sevilla, Spain.}
\footnotetext[7]{Max Planck Institute of Molecular Cell Biology and Genetics, Dresden 01307, Germany.}
\footnotetext[8]{Centre for Systems Biology Dresden, Dresden 01307, Germany.}
\footnotetext[9]{Faculty of Mathematics, Technische Universit\"{a}t Dresden, Dresden 01062, Germany.}
\footnotetext[10]{Corresponding author(s). Maria-Jose Jimenez: majiro@us.es; Heather A. Harrington: harrington@mpi.cbg.de}
\vspace{-3.5\baselineskip}
\begin{abstract}
	Spatial relationships in multi-species data can indicate and affect system outcomes and behaviors, ranging from disease progression in cancer to coral reef resilience in ecology; therefore, quantifying these relationships is an important problem across scientific disciplines.
	Persistent homology (PH), a key mathematical and computational tool in topological data analysis (TDA), provides a multiscale description of the shape of data.
	While it effectively describes spatial organization of species, such as cellular patterns in pathology, it cannot detect the shape relations between different types of species.
	Traditionally, PH analyzes single-species data, which limits the spatial analysis of interactions between different species.
	Leveraging recent developments in TDA and computational geometry, we introduce a scalable approach to quantify higher-order interactions in multi-species data.
	The framework can distinguish the presence of shape features or patterns in the data that are (i) common to multiple species of points, (ii) present in some species but disappear in the presence of other species, (iii) only visible when multiple species are considered together, and (iv) formed by some species and remain visible in the presence of others.
	We demonstrate our approach on two example applications.
	We identify (1) different behavioral regimes in a synthetic tumor micro-environment model, and (2) interspecies spatial interactions that are most significantly altered in colorectal cancer tissue samples during disease progression.
\end{abstract}

Many interactions across the sciences are inherently spatial, involving the relationships between multiple species in concert, yet quantification of such multi-species spatial organization is limited.
Further, the emergent properties of complex systems are often driven by higher-order interactions that cannot be reduced to pairwise connections.
For example, wound healing requires coordination between fibroblasts, endothelial cells, and immune cells for angiogenesis and immune response; and the formation of tertiary lymphoid structures involves the coordination of T cells, B cells, and fibroblasts.
In cancer biology, the spatial organization of and interaction between distinct cell populations within a tumor micro-environment---for example tumor cells, macrophages, and endothelial cells---profoundly affects tumor progression, metastasis, and treatment responses~\cite{anderson2020tumor}.
Developments in multiplex imaging techniques~\cite{goltsev2018deep,pourmaleki2023moving,giesen2014highly} can produce, through the expression of functional markers, spatially resolved maps of tissue samples at single-cell resolution; after post-processing, the result is the coordinates $(x, y)$ of individual cells of each cell type.
Quantifying the spatial interactions among different cell types in these maps can help determine and describe the primary interactions that distinguish healthy and diseased tissue.

A diverse set of approaches can measure multi-species interactions, but at present, none of these offer a practical and interpretable quantification of higher-order topological relationships.
Methods from spatial statistics characterize multi-species interactions in terms of statistical quantities derived from the cells' spatial distributions and their distances~\cite{bull2020combining,hagos2022high,feng2023spatial,Bull2024extended,schladitz_third_2000,kerscher_statistical_2000}.
Network approaches quantify pairwise interactions in terms of connectivity information between the cells~\cite{weeratunga2023single, weeratunga2025temporo}.
Hypergraphs have been used to model multi-way interactions in genomics and metabolomics~\cite{feng_hypergraph_2021,colby_hypernetworks_2024}, but these structures capture combinatorial relationships rather than spatial topological features, and are less well-suited to spatial data.
Topological data analysis (TDA) provides tools to describe multiscale spatial connectivity---like connected components, loops, and voids---the most prominent of which is persistent homology (PH).
By combining the PH from each species, one can study multi-species spatial relationships as in~\cite{vipond2021multiparameter,stolz2022multiscale,yang2025topological}.
More recent TDA methods can analyze pairwise interactions directly, using Dowker persistence and multi-species witness filtrations~\cite{stolz2024relational,yoon2024deciphering}, or mixup barcodes~\cite{wagner2024mixup}.
These topological methods are only applicable to three or more species through the induced pairwise relationships.
Statistical and topological descriptors can also be integrated to leverage their respective strengths, as shown by Bull et al.~\cite{bull2024integrating,Bull2024Muspan} in their analysis of multi-species pairwise spatial interactions in colorectal cancer.

Directly quantifying higher-order interactions in multi-species data remains a significant challenge.
The challenge is both mathematical---to define and interpret higher-order interactions---and computational, since the number of possible interactions grows exponentially with the number of species in the data.
Chromatic TDA, a recent advancement in topological data analysis~\cite{Montesano2024chromaticTDA,Montesano2026chromatic}, addresses many of these issues; however, theoretical and computational gaps remain with the practical application and interpretation of the framework.
Building on this work, we propose a novel topological framework to quantify higher-order spatial topological interactions in multi-species data (summarized in \Cref{fig: pipeline and comparison of methods}).
The framework leverages the discrete Morse theory of chromatic Delaunay filtrations~\cite{Natarajan2024} for computational efficiency, and produces interpretable topological signatures of multi-species spatial relationships.
These signatures can detect the presence of various spatial relationships such as (1) topological features that are present in some species but are masked in the presence of other species (\Cref{fig: example 2a,fig: example 2b}), (2) features that only appear when multiple species are considered together (\Cref{fig: example 2c}), and (3) features that are shared by multiple species (\Cref{fig: example 2d}).
In particular, the method can detect spatial features that are not present when considering only pairwise interactions.

We showcase the proposed framework on synthetic and empirical
spatial biological data (\Cref{sec: description of datasets,sec: results}), where we analyze spatial features arising from singles, doubles, and triples of species of points.
Our results indicate that the framework consistently agrees with existing pairwise analysis, but provides more complete and robust description of the geometry of multi-species spatial interactions.

\section{Topological model}\label{sec: topological model}
\subsection{Preliminaries}\label{sec: preliminaries}
Central to TDA is the idea of interpreting datasets as geometric objects that can be studied with tools from algebraic topology.
In the classical TDA paradigm, we can associate a $d$-dimensional point cloud $X \subset \mathbb{R}^d $ with a \emph{filtered simplicial complex}---a nested sequence of simplicial complexes $F_r(X)$ indexed by a scale parameter $r$---that serves as a multiscale approximation of the underlying space from which $X$ is assumed to be sampled.
At each scale $r$, the $m$-dimensional connectivity of the complex $F_r(X)$---for example connected components, loops, and voids when $m=0, 1$ and $2$ respectively---is described by a vector space $H_m(F_r(X))$ called the \emph{degree-$m$ homology} of $F_r(X)$.
Between any two scale values, say $r \leq s$, the inclusion of complexes $F_r(X) \hookrightarrow F_s(X)$ induces a linear map between the corresponding homology vector spaces $H_m(F_r(X)) \to H_m(F_s(X))$.
The collection of degree-$m$ homology vector spaces across all scale values, and the linear maps between these spaces, together comprise the \emph{degree-$m$ persistent homology} $PH_m(X)$ of $X$~\cite{frosini_measuring_1992,robins_towards_1999,edelsbrunner2000topological,zomorodian2005computing}.
Persistent homology (PH) describes topological features at each scale in the data, and the evolution of these features across scales.
PH can be completely described (and visualized) by its \emph{persistence diagram}, which is a collection of points in the plane where each point $(b, d)$ represents a topological feature which is born at scale $b$ and dies at scale $d$.
The \emph{lifetime} of the feature is $d-b$ (visually, the vertical displacement of $(b, d)$ above the diagonal); long-lived features are typically interpreted as important features of the underlying space, while short-lived features are often interpreted as noise.
Under mild assumptions, $PH_m(X)$ is stable to perturbations in $X$~\cite{cohen2005stability,cohen-steiner_lipschitz_2010}, making it a robust feature descriptor which is suitable as input (after vectorization) for machine learning and statistical inference pipelines.

For multi-species data---wherein $X$ is a disjoint union $\bigsqcup_{i=0}^s X_i$ of points belonging to different species---there have been several attempts to modify the traditional TDA pipeline to incorporate the additional information provided by species labels.
In~\cite{stolz2024relational}, the authors use \emph{Dowker persistence} as defined in~\cite{chowdhury_functorial_2018}, and also introduce \emph{multi-species witness filtrations}.
In this setting, for fixed subspecies of points $X_i$ and $X_j$, a simplicial complex is built from $X_i$ which is filtered by the relative proximity of its simplices to points in $X_j$, and the persistent homology of the filtration is computed.
\emph{Mixup barcodes}, proposed in~\cite{wagner2024mixup}, combine standard persistent homology with image persistent homology (described below) to capture how the lifetimes of topological features in $X_i$ are shortened in the presence of points from $X_j$.
Both of these approaches are restricted to analyzing interactions between pairs of point clouds.

\subsection{Chromatic TDA}
Montesano et al.~\cite{Montesano2026chromatic} propose a framework called \emph{chromatic topological data analysis} that combines kernel, cokernel, and image persistence.
The main idea in this framework is that an inclusion of points $X_i \subseteq X_i \cup X_j$ induces an inclusion of complexes \mbox{$F_r(X_i) \subseteq F_r(X_i \cup X_j)$} at each scale, and these induce linear maps between the homology spaces \mbox{$\iota_r : H_m(F_r(X_i)) \to H_m(F_r(X_i \cup X_j))$}.
These linear maps, collectively across scales, comprise a map of persistent homology modules \mbox{$\iota: PH_m(X_i) \to PH_m(X_i \cup X_j)$}.
The \emph{kernel} of this map---itself a collection of vector spaces $\{\ker(\iota_r)\}_{r \geq 0}$ and linear maps $\{\ker(\iota_r) \to \ker(\iota_s)\}_{r \leq s}$ across scales---describes topological features in $X_i$ that are not present in $X_i \cup X_j$, as well as the evolution of such features across scales.
The \emph{cokernel} of $\iota$ describes features that are present in $X_i \cup X_j$ but not present in $X_i$ alone.
The \emph{image} of $\iota$ describes features in $X_i$ that are still present in $X_i \cup X_j$.
The kernel, cokernel, and image can also be described by persistence diagrams.
Combining these diagrams with the persistent diagrams of $PH_m(X_i)$ (the \emph{domain}) and $PH_m(X_i \cup X_j)$ (the \emph{codomain}), one obtains a rich mathematical description of the spatial relationships between $X_i$ and $X_i \cup X_j$; these are part of the \emph{6-pack} of diagrams proposed in~\cite{Montesano2026chromatic}.
The above approach can be extended to more than two species by replacing $X_i$ and $X_j$ with $\sqcup_{i \in I} X_i$ and $\sqcup_{j \in J} X_j$ for disjoint sets of colors $I$ and $J$.
However, the method is not symmetric in the chosen species $i$ and $j$; this is problematic when $X_i$ and $X_j$ have very different spatial densities.
Moreover, the total number of persistence diagrams increases sharply with the number of species after considering all possible inclusions.

\begin{figure*}
	\begin{subfigure}{\linewidth}
		\centering
		\input{figures/pipeline/pipeline_diagram.tex}
		\caption{
			The \Pipeline pipeline.
			A collection of point clouds with inclusions between them yields a gluing map of filtrations (the $k$-chromatic gluing map), which induces a map of persistent homology.
			This map can be summarized using the domains, kernel, cokernel, and image persistence diagrams.
			Here the domain diagram is the usual persistence diagram of each point cloud separately, while the other diagrams encode the various spatial relationships among the point clouds.
			These diagrams are then vectorized using persistent statistics to obtain feature vectors that summarize spatial relationships among the point clouds.
		}
		\label{fig: pipeline}
	\end{subfigure}\\[0.5\baselineskip]
	\begin{subfigure}[t]{0.4\linewidth}
		\centering
		\includegraphics[width=\linewidth]{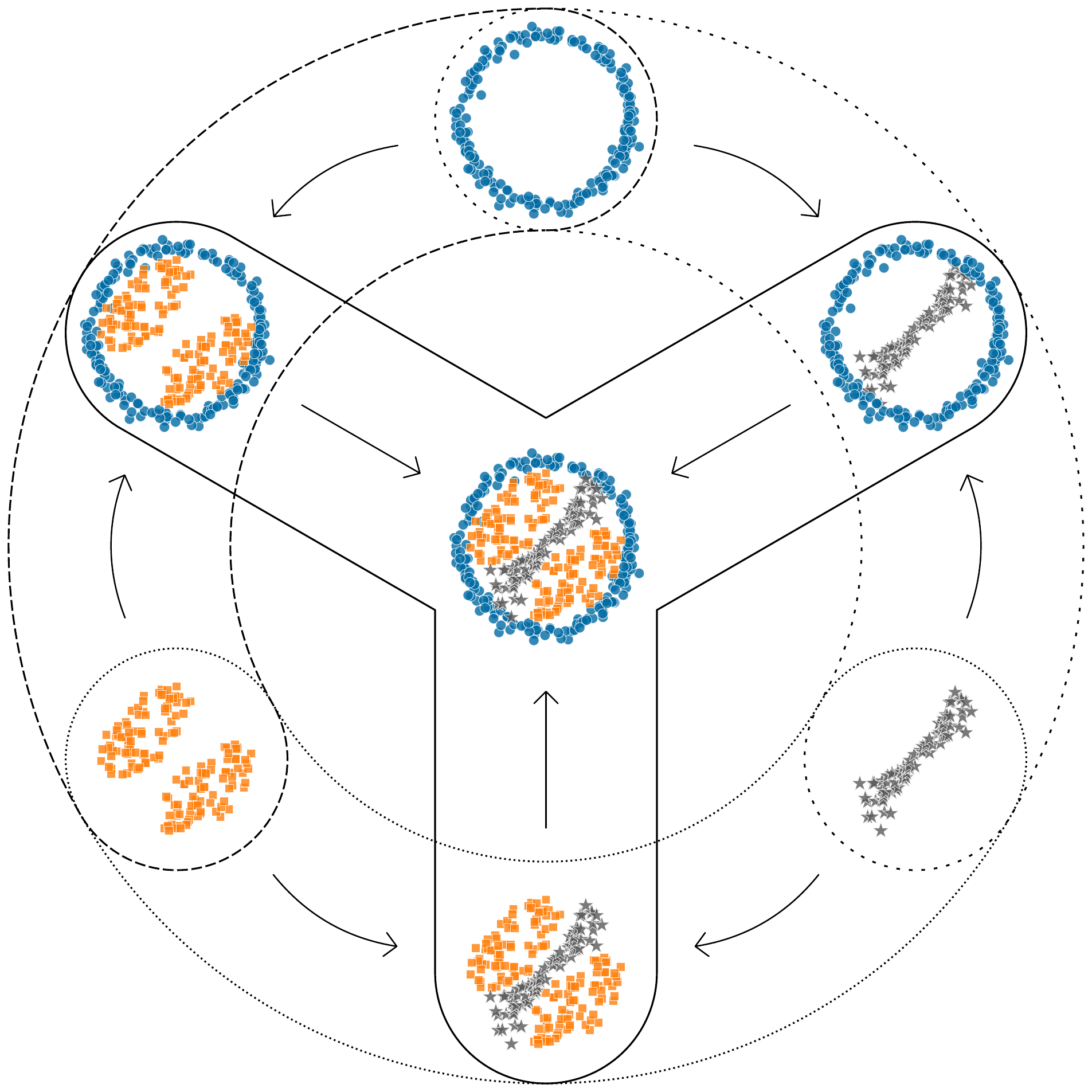}
		\caption{
			The gluing maps used for a trichromatic example.
			The arrows in each dashed/dotted region yield a 1-chromatic gluing map for two species.
			The arrows in the region bounded by the solid line comprise the 2-chromatic gluing map.
			There are three 1-chromatic gluing maps and one 2-chromatic gluing map, yielding four maps in all, each of which is used in the pipeline shown in (a).
		}
		\label{fig: inclusions for three colours}
	\end{subfigure}%
	\hfill%
	\begin{subfigure}[t]{0.57\linewidth}
		\centering
\definecolor{cbBlue}{RGB}{0,107,164} 
\definecolor{Orange}{RGB}{255,128,14} 
\colorlet{cbOrange}{Orange!80} 
\definecolor{cbGrey}{RGB}{89,89,89} 
\newcommand{\cmark}{\ding{51}} 
\newcommand{\xmark}{\ding{55}} 

\begin{tblr}{
	baseline={b},
	width=\linewidth,
	colspec={X[2,c,m]X[c,m]X[c,m]X[c,m]},
	colsep={1pt},
	hlines, 
	hline{1,2,Z} = {1pt}, 
		}
	\small Examples
	 & \small Dowker persistence \cite{stolz2024relational}
	 & \small $k$-chromatic inclusion \cite{Montesano2026chromatic}
	 & \small $k$-chromatic gluing map (proposed)                   \\

	$\vcenter{\hbox{\begin{tikzpicture}[scale=0.5]
					                \draw[black] (-1.4,-1.4) rectangle (1.4,1.4);
					                \draw[line width=2pt, cbBlue] (0,0) circle (1.2);
					                \draw[line width=2pt, cbOrange] (0,0) circle (1);
				                \end{tikzpicture}
				\begin{tikzpicture}[scale=0.5]
					\draw[black] (-1.4,-1.4) rectangle (1.4,1.4);
					\fill[cbOrange] (0,0) circle (1.0);
					\draw[line width=2pt, cbBlue] (0,0) circle (1.2);
				\end{tikzpicture}}}$
	 & \Large \xmark
	 & \smaller Image
	 & \smaller Image                                               \\

	$\vcenter{\hbox{\begin{tikzpicture}[scale=0.5]
					                \draw[black] (-1.4,-1.4) rectangle (1.4,1.4);
					                \draw[line width=2pt, cbBlue] (-0.55,0) circle (0.70);
					                \draw[line width=2pt, cbOrange] (0.85,0) circle (0.50);
				                \end{tikzpicture}
				\begin{tikzpicture}[scale=0.5]
					\draw[black] (-1.4,-1.4) rectangle (1.4,1.4);
					\draw[line width=2pt, cbBlue] (-0.55,0) circle (0.70);
					\draw[line width=2pt, cbOrange] (0.85,0) ++(90:0.50) arc (90:270:0.50);
				\end{tikzpicture}}}$
	 & \Large \xmark
	 & \smaller Image
	 & \smaller Image                                               \\
	$\vcenter{\hbox{\begin{tikzpicture}[scale=0.5]
					                \draw[black] (-1.4,-1.4) rectangle (1.4,1.4);
					                \draw[line width=2pt, cbBlue] (0,0) circle (1.2);
					                \draw[line width=2pt, cbOrange] (0,0) ++(0:1.0) arc (0:180:1.0);
					                \draw[line width=2pt, cbGrey] (0,0) ++(-180:1.0) arc (-180:0:1.0);
				                \end{tikzpicture}
				\begin{tikzpicture}[scale=0.5]
					\draw[black] (-1.4,-1.4) rectangle (1.4,1.4);
					\draw[line width=2pt, cbBlue] (0,0) ++(0:1.2) arc (0:180:1.2);
					\draw[line width=2pt, cbOrange] (0,0) ++(-180:1.2) arc (-180:0:1.2);
					\draw[line width=2pt, cbGrey] (0,0) ++(-180:1.0) arc (-180:0:1.0);
				\end{tikzpicture}}}$
	 & \smaller Needs 3 persistence diagrams
	 & \Large \xmark
	 & \smaller Kernel                                              \\

	$\vcenter{\hbox{\begin{tikzpicture}[scale=0.5]%
					                \draw[black] (-1.4,-1.4) rectangle (1.4,1.4);
					                \draw[line width=2pt, cbBlue]   (-0.55,0) ++(0:0.70) arc (0:180:0.70);
					                \draw[line width=2pt, cbOrange] (-0.55,0) ++(180:0.70) arc (180:360:0.70);
					                \draw[line width=2pt, cbGrey] (0.85,0) circle (0.5);
				                \end{tikzpicture}
				\begin{tikzpicture}[scale=0.5]%
					\draw[black] (-1.4,-1.4) rectangle (1.4,1.4);
					\draw[line width=2pt, cbBlue]   (-0.55,0) ++(0:0.70) arc (0:180:0.70);
					\draw[line width=2pt, cbOrange] (-0.55,0) ++(180:0.70) arc (180:360:0.70);
					\draw[line width=2pt, cbBlue] (0.85,0) ++(0:0.5) arc (0:180:0.5);
					\draw[line width=2pt, cbGrey] (0.85,0) ++(180:0.5) arc (180:360:0.5);
				\end{tikzpicture}}}$
	 & \smaller Needs 3 persistence diagrams
	 & \Large \xmark
	 & \smaller Kernel                                              \\
\end{tblr}%
		\caption{
			Distinguishing power of our approach vis-a-vis other topological methods.
			The first column shows examples of spatial patterns of labelled point clouds.
			Each of the remaining columns corresponds to a topological method to quantify these patterns, indicating whether the method has the ability to distinguish the examples in the first column.
			The second column refers to the use of Dowker persistence diagrams, which cannot deal with trichromatic spatial relationships directly.
			The combination of the 3 pairwise Dowker diagrams can distinguish the examples in the last two rows.
			For the remaining columns we explicitly state which of the diagrams (domain, codomain, kernel, image, cokernel) can distinguish between the configurations, for $k=\#\mbox{ colors} - 1$.
		}
		\label{fig: comparison of methods}
	\end{subfigure}
	\caption{The \Pipeline pipeline (a) combines information from the persistence diagram of each species of points, along with the kernel, cokernel and image diagrams arising from gluing combinations of species together according to the scheme shown in (b).
		This results in a more interpretable and discriminative tool (c).}\label{fig: pipeline and comparison of methods}
\end{figure*}

The authors address these issues by defining the \emph{$k$-chromatic sub-filtration} of $F_{\bullet}(X)$, consisting of all simplices with points of at most $k$ distinct colors.
We refer to the inclusion of this sub-filtration into the full chromatic filtration is the \emph{$k$-chromatic inclusion map}.
The kernel, cokernel, and image diagrams of the induced map on persistent homology can capture many important spatial relationships among the species.
However, some problems remain with this approach.
First, it does not identify the co-occurrence of features in multiple species of points (see \Cref{fig: comparison of methods}); and second, interpreting features in the diagrams associated to this map is not straightforward.
We elaborate on this point in the SI Appendix.

\subsection{Multiscale Multi-species Spatial Signatures (\Pipeline)}\label{sec: maths description}
We propose a framework to address the various shortcomings of existing topological methods for analyzing multi-species data.
For a subset of species $I \subset \{0, \ldots, s\}$, we let $X_I$ denote the subset of points in $X$ whose labels belong to $I$.
Fixing an integer $k$, we consider the disjoint union of the complexes spanned by $X_I$ as $I$ ranges over $k$-ary combinations of colors:
\begin{equation}
	F^{(k)}_r(X) \coloneqq \bigsqcup_{\substack{I \subset \{0, \ldots, s\}\\|I| = k}} F_r(X_I).
\end{equation}
Mapping each simplex in $F_r^{(k)}(X)$ to the same simplex in $F_r(X)$ gives a map \mbox{$F_r^{(k)}(X) \to F_r(X)$}, which in turn yields a map \mbox{$PH_m(F^k(X)) \to PH_m(F(X))$} which we refer to as the \emph{\hbox{$k$-chromatic gluing} in degree $m$}; this map glues together the filtrations spanned by the points in each $k$-ary combination of species.
We consider the kernel, cokernel, and image persistence diagrams of this map.
These diagrams can be vectorized to obtain feature vectors---for example using \emph{persistent statistics}~\cite{ali2023survey} or \emph{persistent landscapes}~\cite{bubenik_statistical_2015}---that describe the spatial interactions among the subsets of the data.

\begin{figure*}[p]
	\centering
	\def\imagewidth{0.48\linewidth}
	\def\captionwidth{0.43\linewidth}
	\begin{subfigure}[t]{\imagewidth}
		\centering
		\includegraphics[valign=t, width=\linewidth]{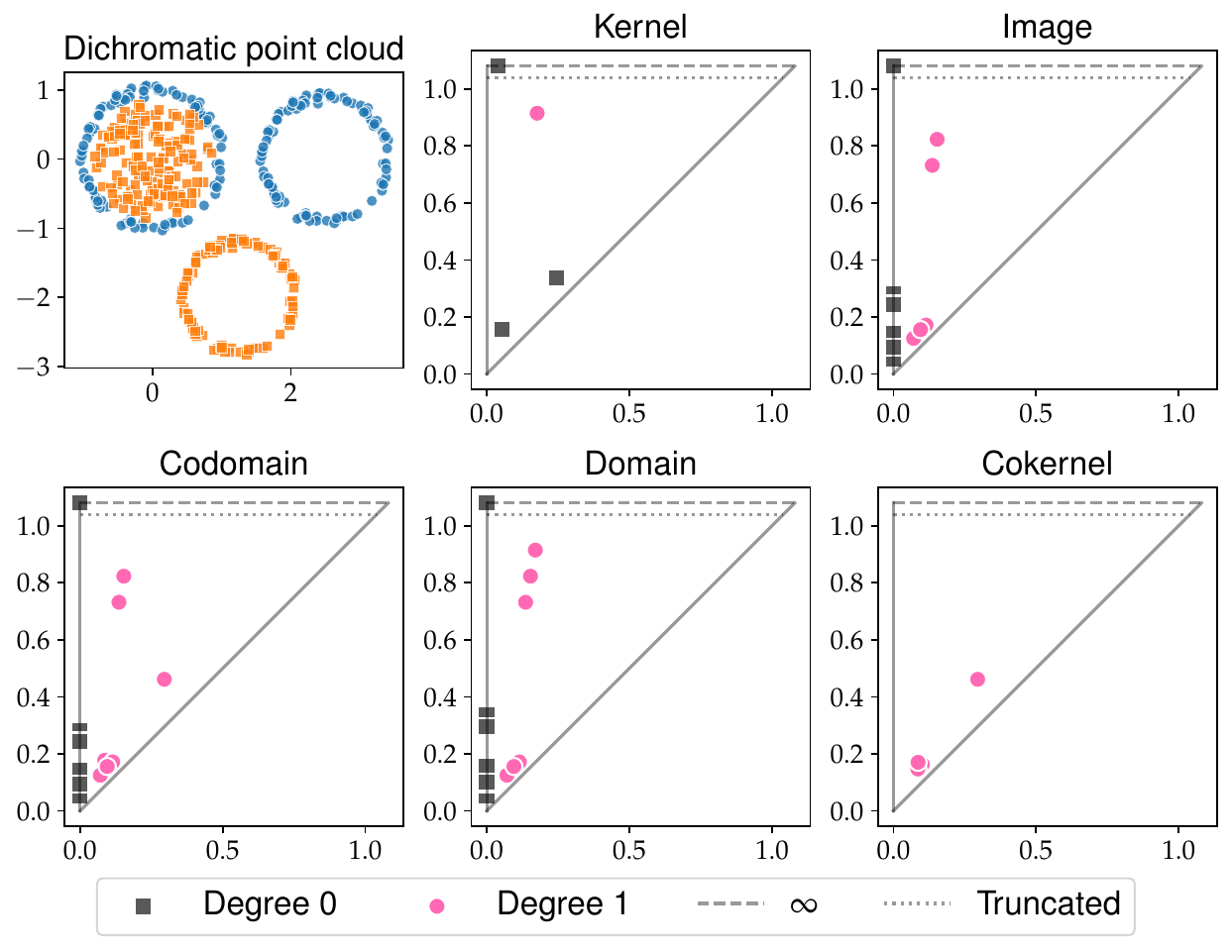}
		\caption{
			A dichromatic point cloud, and the associated kernel, image, codomain, domain and cokernel persistence diagrams arising from the 1-chromatic gluing map.
			The kernel in degree 1 indicates the presence of a monochromatic loop which is filled in by points of the other color.
			The two long-lived features in the image in degree 1 indicate the presence of two monochromatic loops that remain when the points of both colors are considered together.
			The three long-lived features in the domain in degree 1 indicate the presence of three monochromatic loops in total.
			The lack of long-lived features in the cokernel diagram indicates that every dichromatic feature corresponds to some monochromatic feature formed by points in at least one of the colors.
		}
		\label{fig: example 2a}
	\end{subfigure}\hfill
	\begin{subfigure}[t]{\imagewidth}
		\centering
		\includegraphics[valign=t, width=\linewidth]{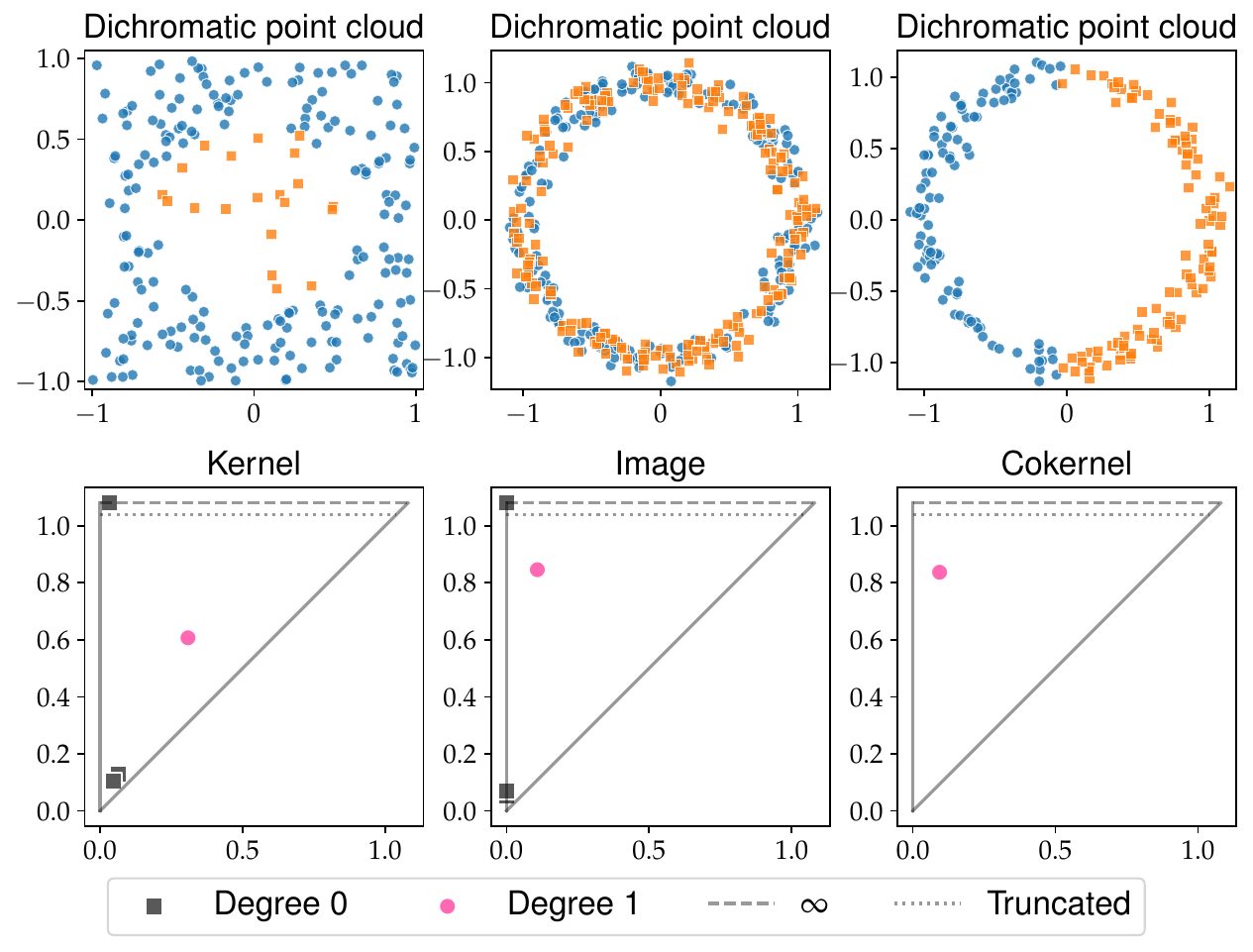}
		\caption{
			When a monochromatic loop is partially or sparsely filled-in by points of a different color, the kernel in degree 1 of the 1-chromatic gluing map captures this phenomenon, as shown in the first example on the left.
			The lifetime of this feature is less than if the loop were completely or densely filled-in; compare with \Cref{fig: example 2a}.
			In the second example, each monochromatic loop remains a loop when the points of the other color are added in, and both monochromatic loops correspond to the same dichromatic loop.
			Therefore, we see a long-lived feature in the image in degree 1.
			In the third example from the left, the singular long-lived feature in the cokernel in degree 1 shows that there is a loop in the dichromatic point cloud that is not present in either of its monochromatic subsets.
		}
		\label{fig: example 2b}
	\end{subfigure}\\[\baselineskip]
	\begin{subfigure}[t]{\imagewidth}
		\centering
		\includegraphics[valign=t,width=\linewidth]{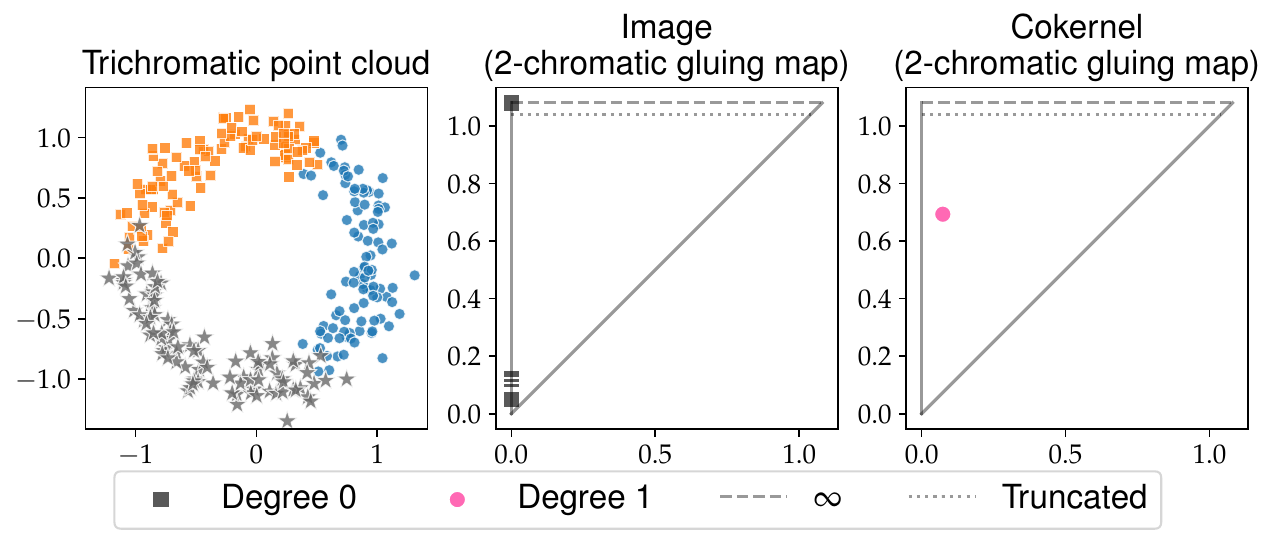}
		\caption{
			For trichromatic point clouds, in addition to diagrams arising from the 1-chromatic gluing maps, we can also consider diagrams arising from the 2-chromatic gluing map.
			In this example, the singular feature in the cokernel in degree 1 attests to the existence of a trichromatic loop that is not present in any dichromatic subset of points.
			The lack of features in the image in degree 1 shows that no large loops are present in any dichromatic subset of points that remain when all three colors are considered together.
			These observations taken together imply that there are no large dichromatic loops in the point cloud.
		}
		\label{fig: example 2c}
	\end{subfigure}\hfill
	\begin{subfigure}[t]{\imagewidth}
		\centering
		\includegraphics[valign=t,width=\linewidth]{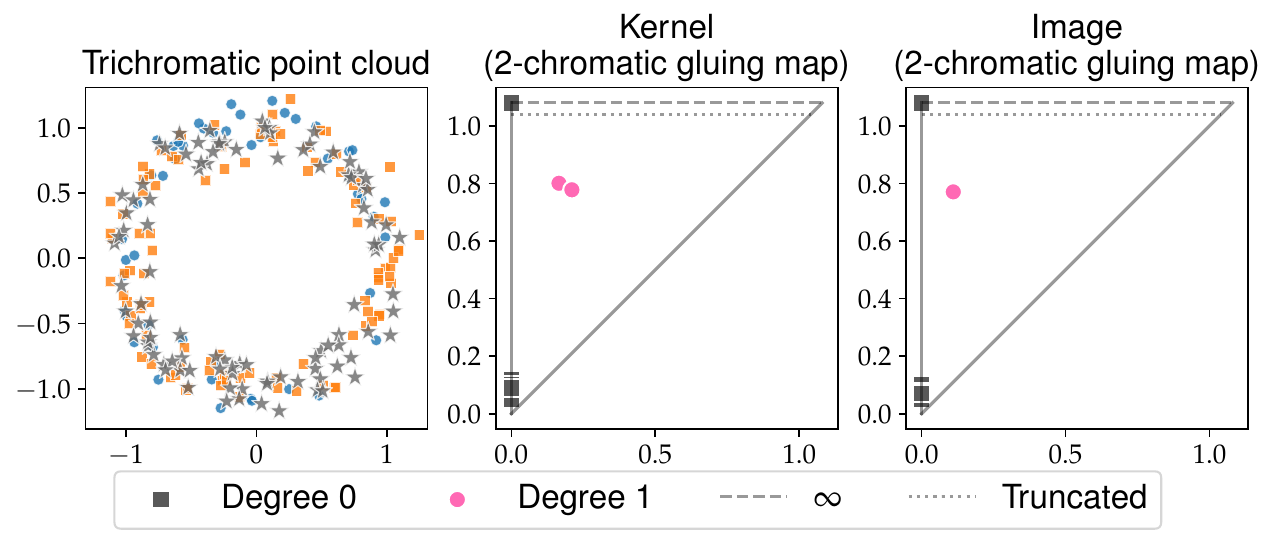}
		\caption{
			In this example, we again consider the 2-chromatic gluing map.
			Each dichromatic subset of points contains a loop, and all three dichromatic loops are co-located and correspond to the same trichromatic loop.
			The two features in the kernel in degree 1 indicate that two of the three dichromatic loops are identified with the third when all three colors are considered together.
			Clearly the trichromatic loop is in the image of the dichromatic loops, as shown by the single long-lived feature in the image in degree 1.
		}
		\label{fig: example 2d}
	\end{subfigure}
	\caption{
		Toy examples illustrating the information captured by the kernel, image, codomain, domain, and cokernel persistence diagrams from the 1-chromatic and 2-chromatic gluing maps.
		Features with persistence less than $0.05$ are omitted from the diagrams for clarity.
	} \label{fig: toy examples}
\end{figure*}

To explain the spatial relationships encoded by these diagrams, we refer to a topological feature as being \emph{$k$-chromatic} if there is a set of $k$ distinct species labels $I$ such that the feature is present in $X_I$.
Then the cokernel of the $k$-chromatic gluing map in degree $m$ identifies \hbox{$m$-dimensional} features that are not $k$-chromatic; that is, these features are not present in any subset of the data consisting of points of at most $k$ colors.
See the third example in \Cref{fig: example 2b} for $k=1$, and \Cref{fig: example 2c} for $k=2$.
The image of the $k$-chromatic map describes $k$-chromatic features that are still present when points of all other species are added in.
See the unfilled circles in \Cref{fig: example 2a} or the second example in \Cref{fig: example 2b}.
The kernel of the $k$-chromatic map identifies $k$-chromatic features that disappear when points of other species are added in.
A feature can disappear if it is filled in by points of another species (in \Cref{fig: example 2a} the blue circle is filled in by orange points), or when two features merge into a single co-located feature shared by different subsets of $k$ species (the two co-located monochromatic circles in \Cref{fig: example 2b} merge, and the three co-located dichromatic circles in \Cref{fig: example 2d} merge into a single trichromatic circle).
We refer to \Cref{fig: toy examples} for more detailed examples.

An important technical detail is the construction of filtrations from the data points.
Traditional choices in TDA are the \v{C}ech and Vietoris--Rips filtrations, which are unsuitable in the present context due to the large number of simplices in these filtrations, or the alpha filtration~\cite{Edelsbrunner1993alpha}, which has fewer simplices but is also unsuitable here since an inclusion of point clouds does not necessarily induce an inclusion of alpha filtrations.
The \emph{chromatic alpha filtration}~\cite{Montesano2026chromatic} was proposed to address the latter shortcoming.
We instead use the \emph{chromatic Delaunay--\v{C}ech filtration} which is topologically equivalent to the chromatic alpha filtration and more computationally efficient~\cite{Natarajan2024}.

\section{Description of datasets}\label{sec: description of datasets}
\subsection{Synthetic data from tumor micro-environment model}\label{sec: ABM data description}
The Agent-Based Model (ABM) introduced in~\cite{bull2023quantification} simulates the spatial and temporal dynamics of a tumor micro-environment through interactions among tumor, stromal, and necrotic cells, blood vessels, and macrophages.
Macrophages are further divided into antitumor ($M_1$) and pro-tumor ($M_2$) states.
The simulations explore the effects of varying two parameters regulating macrophage behavior: $\chi^m_c$, the chemotactic sensitivity to CSF-1 gradients, and $c_{1/2}$, the CSF-1 concentration at which macrophage extravasation is half-maximal.
Distinct parameter regimes are identified in which the ABM exhibits one of three qualitative behaviors corresponding to the classical phases of cancer immuno-editing~\cite{Dunn2004}: \emph{elimination}, where $M_1$ macrophages eradicate the tumor; \emph{equilibrium}, where macrophages contain but do not destroy the tumor; and \emph{escape}, where $M_2$ macrophages promote tumor migration towards vasculature.
We consider nine distinct values for each parameter, yielding $81$ parameter combinations, with up to $20$ stochastic realizations per pair, yielding $1485$ simulated point clouds generated at the final simulation time ($t = 500$ hours).
Following~\cite{bull2023quantification}, the behavioral outcome for each choice of parameters was qualitatively assigned based on the spatial distributions of the different cell types in the corresponding simulations; see \Cref{fig: ABM examples} for some examples.

\subsection{Empirical data from colorectal cancer specimens}\label{sec: empirical dataset description}
The data comprises $43$ carcinoma-in-adenoma specimens from the \emph{Oxford Rectal Cancer Cohort}, representing colorectal tumors transitioning from adenoma to carcinoma.
Multiplex immunofluorescence staining and image segmentation were used to identify the coordinates of eleven species of points in each specimen based on marker expression~\cite{bull2024integrating}: epithelial cells, neutrophils, macrophages, cytotoxic T-cells, T-helper cells, T-regulatory cells, periostin, CD146, CD34, SMA, and podoplanin.
Each specimen comprises multiple regions of interest (ROIs), classified as either adenoma (benign) or carcinoma (malignant) based on histopathological assessment.
After quality control, 12030 well-aligned ROIs were obtained from 41 of the patient specimens, with 7360 carcinoma ROIs and 4670 adenoma ROIs, ranging from 109 to 700 ROIs per patient (inclusive).
Each ROI is represented as a labelled point cloud in $\mathbb{R}^2$.

\section{Results}\label{sec: results}
We showcase the effectiveness of Multiscale Multi-species Spatial Signatures (\Pipeline) on point locations from two datasets: (1)~data from a computational model of the tumor microenvironment~\cite{bull2023quantification} with five types of cells, and (2)~empirical data from colorectal cancer specimens with eleven species.
In the synthetic dataset, multi-species spatial relationships can distinguish between different qualitative model behaviors.
In the colorectal cancer dataset, multi-species spatial relationships can identify combinations of species whose spatial relationships are altered during the progression from adenoma to carcinoma.
These include three-way spatial relationships not captured by pairwise relationships alone.

\subsection{\texorpdfstring{\Pipeline}{\PipelineText} distinguishes qualitative behavioral regimes in the synthetic dataset}

\begin{figure*}[t]
	\centering
	\begin{subfigure}[t]{\linewidth}
		\centering
		\includegraphics[width=0.7\linewidth]{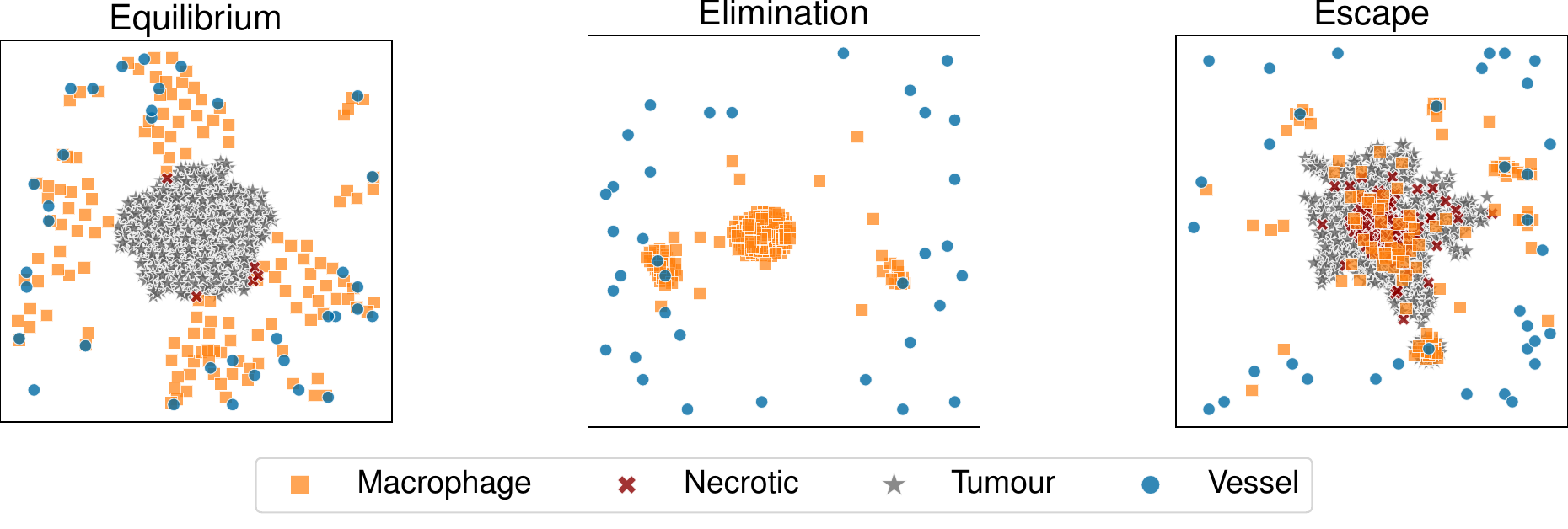}
		\caption{
			Example samples from each behavioral regime of the ABM model.
			Stromal cells are omitted for clarity; M1 and M2 macrophages are not distinguished in the figure and in our analysis.
		}
		\label{fig: ABM examples}
	\end{subfigure}\\[2\baselineskip]
	\begin{subfigure}[t]{0.4\linewidth}
		\centering
		\includegraphics[height=4cm]{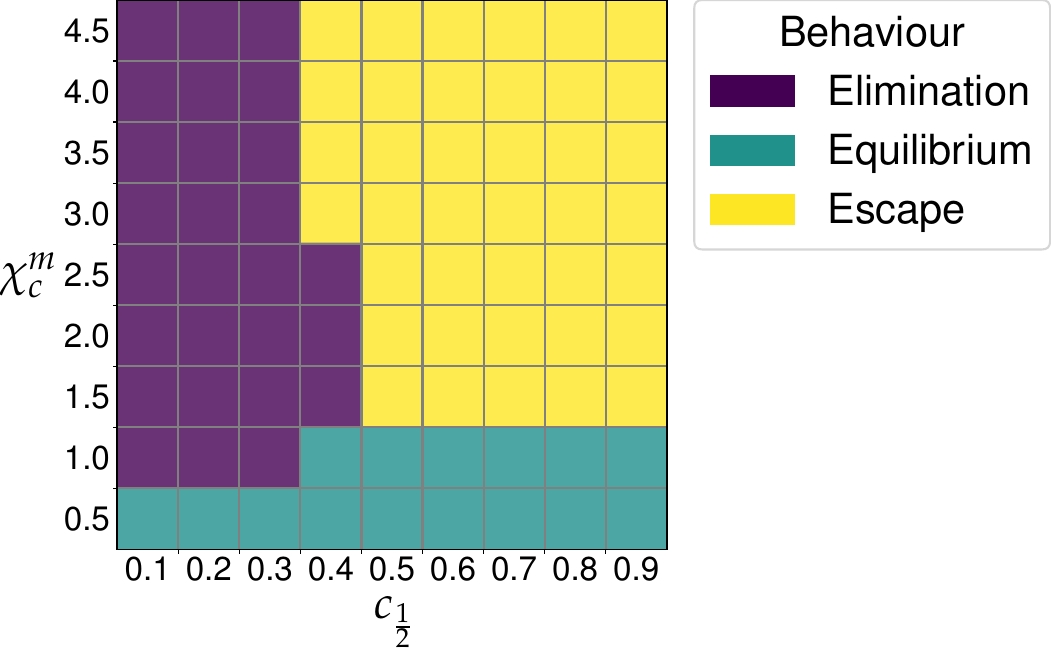}
		\caption{Subjective classification of the ABM parameters.}
		\label{fig: ABM subjective labels}
	\end{subfigure}\hfill
	\begin{subfigure}[t]{0.3\linewidth}
		\centering
		\includegraphics[height=4cm]{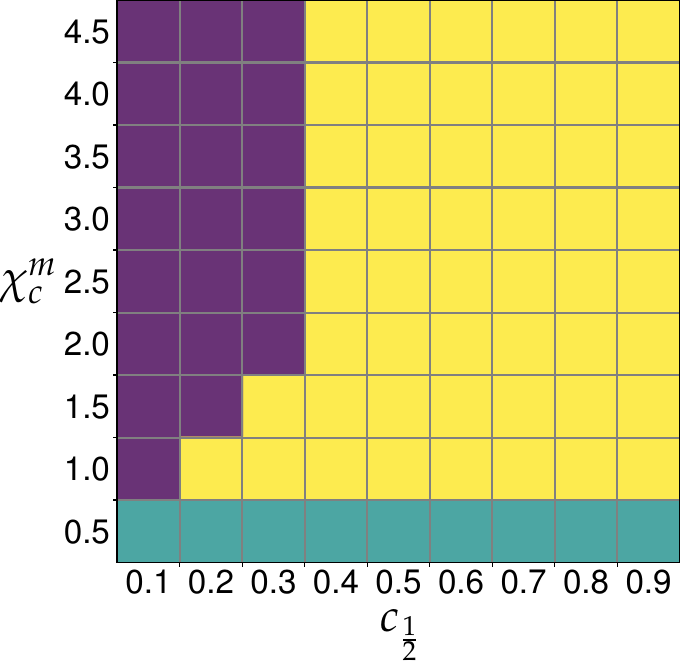}
		\caption{Clustering results using $k$-means with $k=3$.}
		\label{fig: ABM clustering results}
	\end{subfigure}\hfill
	\begin{subfigure}[t]{0.3\linewidth}
		\centering
		\includegraphics[height=4cm]{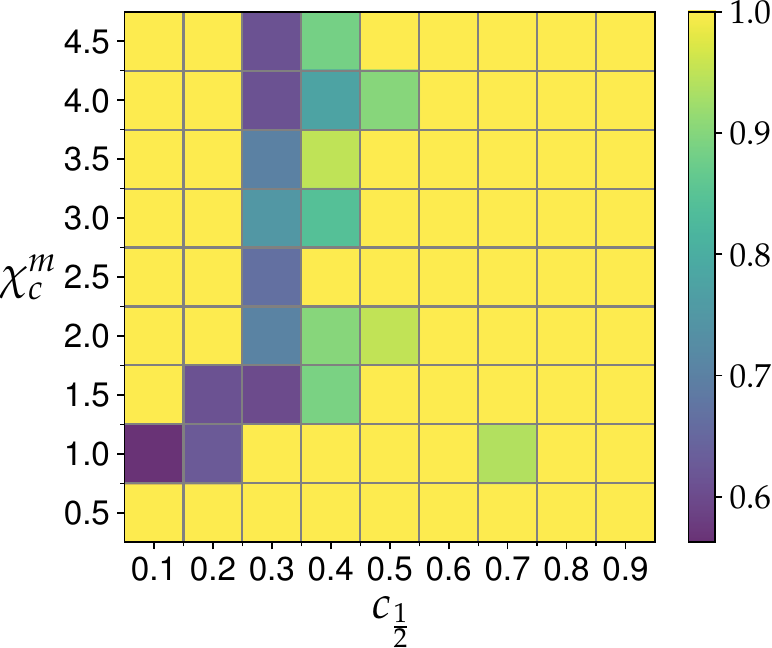}
		\caption{Purity scores of the clustering.}
		\label{fig: ABM purity scores}
	\end{subfigure}
	\caption{
		Clustering results and clustering purity scores for the synthetic ABM dataset.
		Simulations across all combinations of the ABM model parameters $c_{\frac{1}{2}}$ and $\chi_c^m$ were clustered based on the spatial signatures obtained from our pipeline at the final simulation time ($t=500$).
		The color for each parameter combination represents the cluster label assigned to the majority of simulations from that parameter combination.
		Purity scores are computed as the proportion of simulations in the majority class for each parameter combination.
	}
	\label{fig:ABM}
\end{figure*}

We use the spatial relationships among the five cell types in each simulation (combining M1 and M2 macrophages) to identify combinations of model parameters that lead to each of the three distinct behavioral outcomes---elimination, equilibrium, and escape (\Cref{fig: ABM examples}).
More precisely, from each simulation we obtain a multi-species spatial signature, and these signatures are assigned to three clusters using unsupervised machine learning.
Each combination of model parameters is labelled by the most frequently occurring cluster label given to the set of simulations generated from that parameter combination.
The resulting clustering (\Cref{fig: ABM clustering results}) is consistent with subjective assignments (\Cref{fig: ABM subjective labels}).
We also examine the clustering purity scores (\Cref{fig: ABM purity scores}) to study the behavior of the model in the transition between behavioral regimes.
While there is a clear boundary between elimination and escape, no boundary region is visible between equilibrium and escape, suggesting that this change in model behavior is more abrupt than the transition from elimination to escape.

\subsection{\texorpdfstring{\Pipeline}{\PipelineText} identifies important spatial interactions in the colorectal cancer dataset}
\begin{figure*}
	\centering
	\begin{subfigure}{\linewidth}
		\centering
		\includegraphics[width=0.75\linewidth]{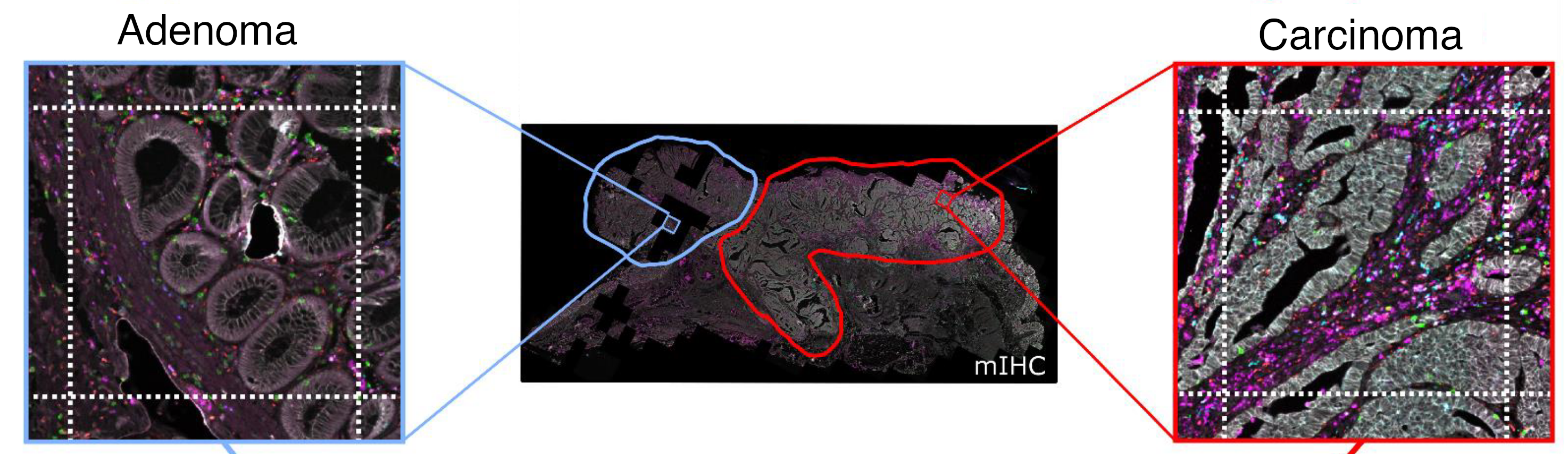}
		\caption{
			Example ROIs from a carcinoma-in-adenoma tissue sample (middle) in the colorectal cancer dataset.
		}
		\label{fig: colorectal cancer example ROI}
	\end{subfigure}\\[\baselineskip]
	\begin{subfigure}{\linewidth}
		\centering
		\includegraphics[width=0.85\linewidth]{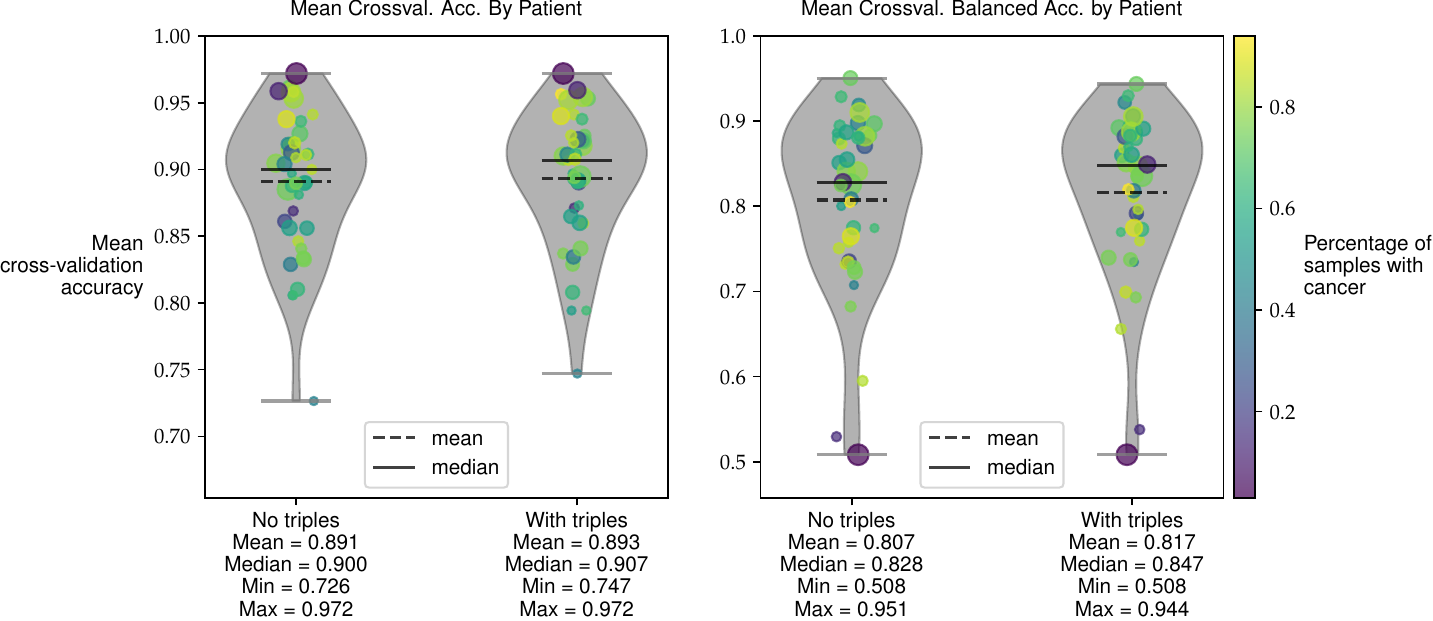}
		\caption{
			The mean classification accuracy across cross-validation folds for each patient in the colorectal cancer dataset, with and without including three-way interactions.
			Each point represents a single patient, and the shaded region is the violin plot of the distribution of accuracies across all patients.
			The area of each point indicates the number of ROIs available for that patient, and the color indicates the number of cancerous ROIs for that patient.
			The left plot shows raw classification accuracies, while the right plot shows the balanced accuracies (the average of the true positive rate and true negative rate) to account for class imbalance.
		}
		\label{fig: colorectal cancer classification accuracy}
	\end{subfigure}\\[\baselineskip]
	\begin{subfigure}{\linewidth}
		\centering
		\includegraphics[height=7cm]{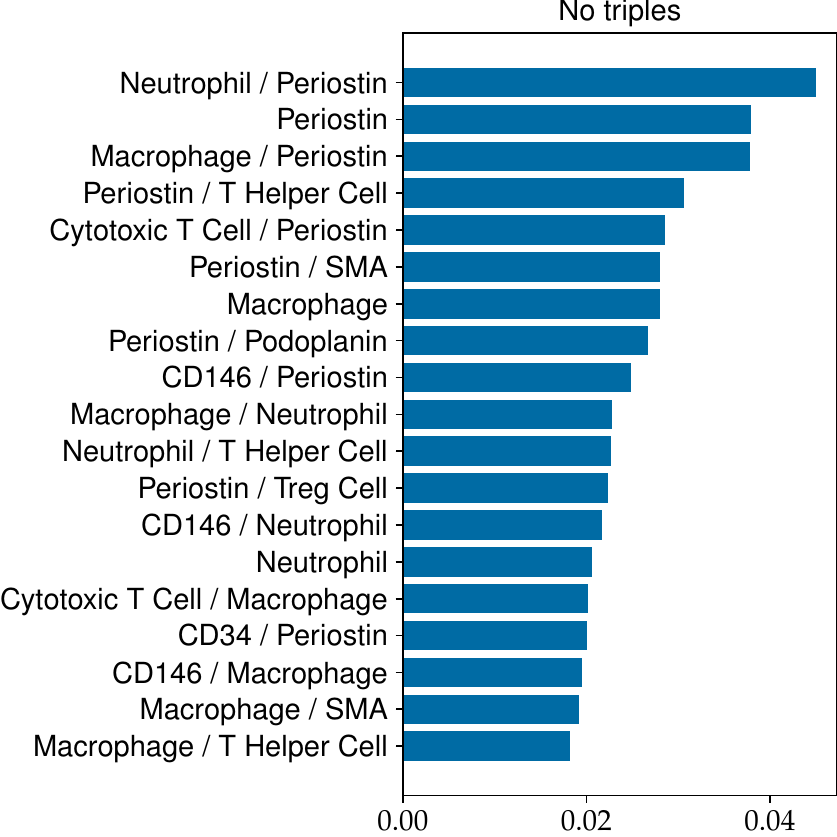}
		\quad
		\includegraphics[height=7cm]{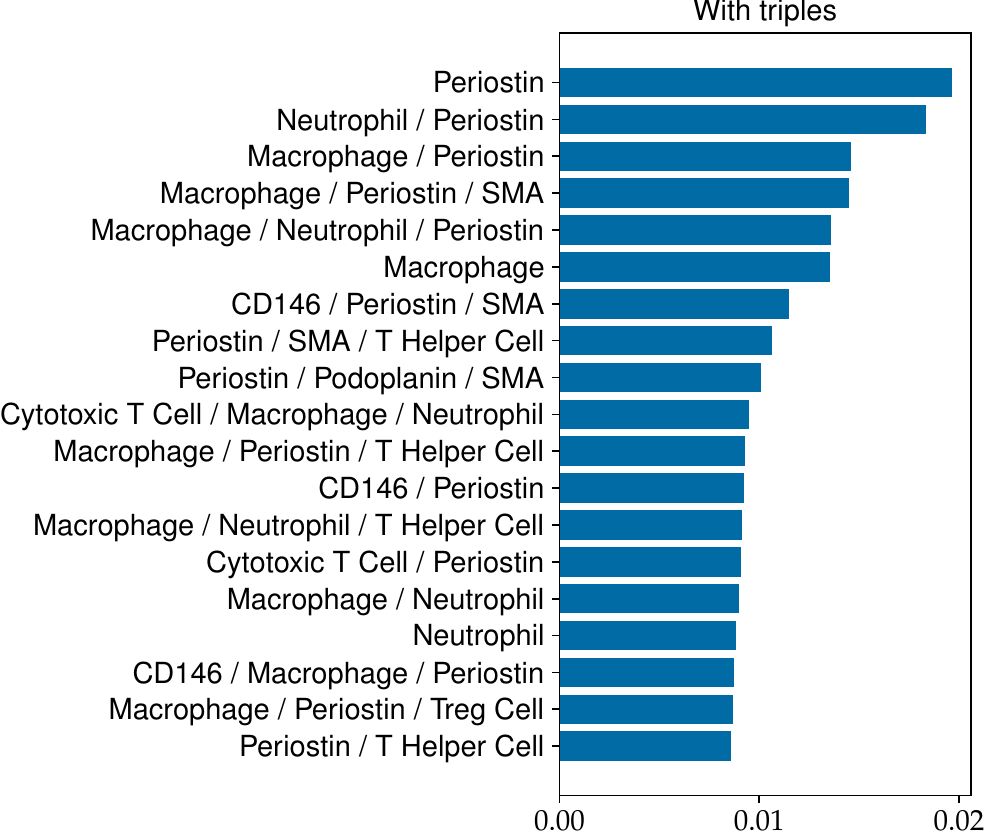}
		\caption{
			The cumulative feature importance scores for each group of species, averaged across all patients in the colorectal cancer dataset, without including three-way interactions (left) and including three-way interactions (right).
			Species groups are sorted by feature importance in descending order, and the top twenty groups are shown in each sub-figure.
			These feature importances are computed using the Mean Decrease in Impurity (MDI) from random forest classifiers trained to distinguish adenoma from carcinoma ROIs for each patient.
		}
		\label{fig: colorectal cancer feature importances}
	\end{subfigure}
	\caption{
		Results of the classification of adenoma and carcinoma ROIs in the colorectal cancer dataset.
	}\label{fig: colorectal cancer results}
\end{figure*}

We apply our pipeline to ROIs obtained from a collection of carcinoma-in-adenoma tissue samples with 11 species of points (\Cref{fig: colorectal cancer example ROI}).
The goal is to identify combinations of species whose spatial relationships are altered during the progression from adenoma to carcinoma, and to assess the importance of three-way spatial relationships over pairwise relationships alone.
For each ROI, we compute multi-species spatial signatures for the 5 immune and 5 stromal species (i.e., excluding the widespread fiducial epithelial marker used for alignment; see \cite{Bull2024extended}), considering combinations of up to three species at a time.
In order to quantify how spatial relationships are altered between adenoma and carcinoma, we train a classifier to distinguish adenoma from carcinoma ROIs based on their multi-species spatial signatures, and we examine the associated feature importance scores.
We find that a classifier trained on the ROIs from all patients together performed worse than a classifier trained separately for each patient due to inter-patient heterogeneity in the spatial structure of the ROIs; therefore, we train a separate classifier for each patient.
We use these classifiers to compute aggregate importance scores for every combination of species considered, averaged across all patients.
The analysis is performed twice: once ignoring triples of species, and once including them.
We emphasize that the goal of our analysis is not the classification of adenoma and carcinoma ROIs per se, but rather the identification of important spatial interactions that are altered during the progression from adenoma to carcinoma.

The prominent periostin--macrophage pairwise spatial interaction uncovered by the analysis is consistent with prior findings in the literature, where periostin is thought to play a role in recruiting tumor-associated macrophages (TAMs) to the tumor site~\cite{zhou_periostin_2015,tang_cross-talk_2018,wei_periostin_2023}.
Including three-way interactions during classification does not significantly affect the classification accuracy (\Cref{fig: colorectal cancer classification accuracy}), but it does uncover potentially important three-way spatial interactions, and it alters the feature importance scores of pairs and singles (\Cref{fig: colorectal cancer feature importances}).
We note in particular the high scores of the macrophage--periostin--neutrophil and macrophage--periostin--SMA three-way interactions.

\section*{Discussion}
We propose a topological data analysis framework to quantify higher-order spatial relationships involving multiple species of points in multi-species spatial data.
The framework distinguishes the presence of topological features in the data---such as connected components, loops, and voids---and the combinatorial arrangement of these features across different subsets of species.
The method leverages recent developments in chromatic topological data analysis and computational geometry, combining kernel, cokernel, and image persistence with a new construction called the \emph{$k$-chromatic gluing map}.
The framework is suitable for analysing higher-order spatial interactions involving three or more species, without reducing these to a set of pairwise interactions as is commonly the case in existing literature.

We validate our framework on synthetic data from an agent-based model of the tumor micro-environment, and empirical data from colorectal cancer tissue samples.
Our results on the synthetic dataset are consistent with existing analysis by Stolz et al.~\cite{stolz2024relational}.
In the colorectal cancer dataset, our results reveal (1) the importance of periostin--macrophage interactions,  which is consistent with existing literature
~\cite{zhou_periostin_2015,tang_cross-talk_2018,wei_periostin_2023}, as well as (2) potentially important three-way interactions involving periostin, macrophages, neutrophils and SMA.
Most importantly, our method is extensible to higher-order interactions and not limited to only pairwise interactions.
We anticipate that our method will be useful in a variety of biological contexts where multi-species spatial data is abundant, such as in the study of microbial communities~\cite{faust_microbial_2012}, tissue architecture \cite{torres-cano_spatially_nodate}, and ecological systems \cite{tarnita_self-organization_2024}.

\section{Materials and Methods}
\subsection{Persistent homology computations}
For each combination of species under consideration, we compute the chromatic Delaunay--\v{C}ech filtration~\cite{Natarajan2024} of the points belonging to those species, and compute persistent homology from the filtration.
For each individual species, we compute 1-parameter persistent homology in degrees 0 and 1.
With pairs (triples, respectively) of species, we consider the 1-chromatic (2-chromatic, respectively) gluing map described in \Cref{sec: maths description}, and we compute the kernel (degrees 0 and 1), image (degrees 0 and 1) and cokernel (degree 1) persistence diagrams associated to this map.
We omit degree 0 cokernel diagrams, since these are always trivial for point clouds.
We do not compute persistent diagrams in degrees greater than 1, since these are always empty for two-dimensional point clouds.

\subsection{Vectorization of persistence diagrams}
We vectorize persistence diagrams using \emph{persistent statistics}~\cite{ali2023survey}, which are statistical summaries built from the birth times, death times, lifespans (death time minus birth time) and midlife points (average of birth and death times) of the features in each persistence diagram.
For each persistence diagram we include the mean, standard deviation, median, range, and the 10th, 25th, 75th and 90th percentiles of birth times, death times, lifespan and midlife point of the features in that diagram.
We also include the number of features and the persistent entropy~\cite{Chintakunta2015pers_entropy1, Rucco2016pers_entropy2} of the diagram.
Thus, we obtain a 34-dimensional feature vector from the persistence diagram.
We exclude the statistics of birth times, midlife and lifetimes for the domain and image diagrams in degree 0, since connected components are always born at time zero, so their lifetime and midlife can be inferred from their death times.
For these diagrams we obtain a 10-dimensional feature vector.
\Cref{table: feature vector dimensions} summarizes the contribution of features from each combination of species.

\begin{table}[H]
	\centering
	\begin{tabular}{|c c c c|}
		\hline
		\textbf{Combination} & \makecell{\textbf{Number of}                                                 \\\textbf{diagrams}} & \makecell{\textbf{Diagram}\\\textbf{types}} & \makecell{\textbf{Feature vector}\\\textbf{dimension}} \\
		\hline
		Individual           & 2                            & Domain (deg. 0,1)            & $10 + 34 = 44$ \\
		\hline
		Pair                 & 5                            & \makecell{Kernel (deg. 0, 1)                  \\Image (deg. 0, 1)\\Cokernel (deg. 1)} & 10 + 4 * 34 = 146 \\
		\hline
		Triple               & 5                            & \makecell{Kernel (deg. 0, 1)                  \\Image (deg. 0, 1)\\Cokernel (deg. 1)} & 10 + 4 * 34 = 146\\
		\hline
	\end{tabular}
	\captionof{table}{Dimension of vectorization of the persistence diagrams.}\label{table: feature vector dimensions}
\end{table}

\subsection{ABM simulations}
The model is a two-dimensional, off-lattice, hybrid, force-based ABM in which cell movement arises from mechanical interactions and chemotactic responses to diffusible species (oxygen, CSF-1, TGF-$\beta$, CXCL12, and EGF)~\cite{bull2023quantification}.
Blood vessels are modelled as stationary and do not interact physically with other cells, while all remaining cell types are motile and interact according to local environmental signals.
Macrophages are characterized by a continuous phenotype \mbox{$\Omega \in [0,1]$}, interpolating between antitumor ($\Omega \approx 0$) and pro-tumor ($\Omega \approx 1$) states.

\subsection{Clustering of ABM simulations}
The input to our pipeline comprises 5 individual cell types, 10 pairs of cell types, and 10 triplets of cell types, for a total of 25 combinations of cell types per simulation.
The feature vectors from these combinations are concatenated to yield a single 3140-dimensional feature vector for each simulation.
These feature vectors are used as input to a $k$-means clustering algorithm with $k=3$ and at most 100000 iterations.

\subsection{Empirical dataset acquisition}
43 carcinoma-in-adenoma specimens were identified from the Oxford Rectal Cancer Cohort and FFPE blocks retrieved following ethical approval (REC 11/YH/0020).
Multiplex immunofluorescence (mIF) staining was conducted using the Akoya Biosciences OPAL Protocol (Marlborough, MA) and the Leica BOND RXm auto-stainer (Leica Microsystems, Germany). Two serial mIF sections per sample were cut and stained to identify epithelial (E-Cadherin), immune (CD4, CD8, FoxP3, CD68, MPO), and mesenchymal (CD34, CD146, Periostin, Podoplanin, $\alpha$-SMA) markers---chosen based on literature for their prognostic relevance and focus on immune–stromal interactions.
The images thus obtained were post-processed through cell segmentation (HALO, Indica Labs), and species were assigned based on marker expression.
Post-processed images were divided into \mbox{1000px $\times$ 1000px} (500$\mu$m$\times$500$\mu$m) regions of interest (ROIs) and aligned using epithelial E-Cadherin as a reference marker present on both sections using a modified version of the \underline{D}ensity \underline{A}daptive Point Set \underline{Re}gistration (DARE) algorithm~\cite{lawin_density_2018}.
See \cite{bull2024integrating} for further details.

\subsection{Classification of adenoma and carcinoma in empirical dataset}
The data comprises 10 individual species (excluding epithelium), 45 pairs of species, and 120 triplets of species, for a total of 175 combinations of species per sample.
The feature vectors from these combinations are concatenated to yield a single $24530$-dimensional feature vector for each ROI.
Due to the large number of species combinations, these feature vectors are very high-dimensional relative to the number of available samples per patient.
Therefore, before training the classifier, we reduce the dimensionality of the feature vectors using Principal Component Analysis (PCA).
For each combination of species $\mathcal{C}$, letting $n_{\mathcal{C}}$ denote the dimension of the feature space corresponding to the features from $\mathcal{C}$, we consider $\min(5, \lfloor 0.05 n_{\mathcal{C}} \rfloor)$ principal components (computed separately for each patient).
This results in 2535 principal components across all the selected combinations of species.
The projection of the data onto these principal components is used as input for the classifier, which is a random forest with 500 gradient-boosted decision trees~\cite{breiman_random_2001,friedman_greedy_2001}.
The maximum tree depth, maximum number of leaf nodes, and minimum number of samples per leaf node were set to 3, 6, and 5 respectively to reduce overfitting, and 3\% of features were sampled for each split.
Gradient boosting was regularized by shrinkage, using a learning rate of 0.4, with 25 iterations of boosting per tree.
The classifier was trained and evaluated using stratified 5-fold cross-validation, with a 70-30 train-test split.
Feature importances are computed using the Mean Decrease in Impurity (MDI)~\cite{breiman_random_2001}.
The net importance score assigned to each group of species is the sum of the MDI scores over all features corresponding to that group, averaged across all patients in the cohort.

\subsection{Dealing with missing species}
In both the synthetic and empirical datasets, for each simulation/sample, we discard species with fewer than 3 representatives to avoid technical difficulties in computing the chromatic Delaunay triangulation.
Features corresponding to species combinations that are not present in a simulation/sample are set to zero in the concatenated feature vector for that simulation/sample so that the feature vectors have a consistent dimension across all simulations/samples.

\subsection{Code}
For computing chromatic filtrations and plotting persistence diagrams, we use version 14.0.0 of the \texttt{chalc} Python package (PyPI project page: \url{https://www.pypi.org/project/chalc}, Github source code: \url{https://github.com/abhinavnatarajan/chalc}).
Persistent homology computations in \texttt{chalc} are performed using version 0.11.0 of the \texttt{lophat} Rust crate (Project page: \url{https://crates.io/crates/lophat}, GitHub source code: \url{https://github.com/tomchaplin/lophat}), which implements the lock-free algorithm for computing persistent homology~\cite{morozov_towards_2020} with the clearing optimization~\cite{bauer_clear_2013}.
Kernel, cokernel, and image diagrams are computed with version 0.3.0 of the \texttt{phimaker} Python package (PyPI project page: \url{https://www.pypi.org/project/phimaker}, GitHub source code: \url{https://github.com/tomchaplin/phimaker}), which uses the algorithm described in~\cite{Cohen-Steiner2009kernels}.
Random forest classification and $k$-means clustering were performed using v1.6 of the Python package \texttt{scikit-learn}~\cite{pedregosa_scikit-learn_2011} (DOI: \href{https://doi.org/10.5281/zenodo.14340513}{10.5281/zenodo.14340513}).

Code for computing the feature vectors from the datasets, running the classifiers, and visualizing the results is available at \url{https://github.com/abhinavnatarajan/\PipelineText_demo}.

\subsection*{Data Archival}
ABM simulations can be found at \url{https://github.com/JABull1066/MacrophageSensitivityABM/}.
Point clouds from the carcinoma-in-adenoma dataset are available at \url{https://github.com/JABull1066/SHIFT-Score-Carcinoma-in-Adenoma/}.

\subsection*{Supporting Information Appendix (SI)}
The SI Appendix contains mathematical background on topological data analysis, and additional figures (Figures S1--S5) to illustrate the mathematical framework:
\begin{enumerate}
	\item Figure S1 shows an example to highlight the difference between our framework and \emph{$k$-chromatic subfiltrations}.
	\item Figure S2 compares our framework with Dowker persistence, expanding on the examples in the first two rows of the table in \Cref{fig: comparison of methods}.
	\item Figures S3 and S4 expand on the examples in the third and fourth rows of the table in \Cref{fig: comparison of methods}, comparing the persistence diagrams obtained from the $k$-chromatic gluing map with those from the $k$-chromatic subfiltrations.
	\item Figure S5 expands upon the examples from \Cref{fig: toy examples}, showing the remaining persistence diagrams for some of those examples.
\end{enumerate}

\subsection*{Acknowledgements}
The authors thank U. Tillmann for helpful discussions and comments on the manuscript.
A.N. is supported by the Clarendon Fellowship, Merton College and the Mathematical Institute at the University of Oxford.
J.A.B. is supported by Cancer Research UK (CR-UK) grant number CTRQQR-2021\textbackslash100002, through the Cancer Research UK Oxford Centre.
E.J.M.I is funded by the Lee Placito Medical Research Felloship (University of Oxford).
S.J.L. is supported by CRUK Programme Grant (DRCNPG-Jun22/100002).
Mathematical and biological integration is supported by CRUK CRC-STARS Strategic Grant (SEBCRCS-2024/100001).
M.J.J. is supported by grant PID2024-155867NB-I00 funded by MICIU/AEI/10.13039/501100011033 and ERDF/EU and by the IMUS–María de Maeztu grant CEX2024-001517-M (Apoyo a Unidades de Excelencia María de Maeztu), funded by MICIU/AEI/10.13039/501100011033.
A.N., T.C., H.M.B., M.J.J. and H.A.H. are grateful for the support provided by the UK Centre for Topological Data Analysis EPSRC grant EP/R018472/1.
HAH gratefully acknowledges funding from the Royal Society RGF/EA/201074, UF150238 and EPSRC EP/Y028872/1 and EP/Z531224/1.

\subsection*{Declaration of Interests}
The authors declare no competing interests.

\printbibliography

\end{document}


\maketitle
\vspace{-2\baselineskip}
\footnotetext[1]{Mathematical Institute, University of Oxford, Radcliffe Observatory Quarter, Andrew Wiles Building, Woodstock Rd, Oxford OX2 6GG, UK.}
\footnotetext[2]{Wolfson Centre for Mathematical Biology, Mathematical Institute, University of Oxford, Oxford OX2 6GG, UK.}
\footnotetext[3]{Centre for Human Genetics, Roosevelt Drive, University of Oxford, Oxford, UK.}
\footnotetext[4]{Ludwig Institute for Cancer Research, Nuffield Department of Medicine, University of Oxford, Oxford OX3 7DQ, UK.}
\footnotetext[5]{Translational Gastroenterology Unit, John Radcliffe Hospital, University of Oxford, and Oxford NIHR Biomedical Research Centre, Oxford, UK.}
\footnotetext[6]{Universidad de Sevilla. Escuela Técnica Superior de Ingeniería Informática. Av. Reina Mercedes s/n, 41012 Sevilla, Spain.}
\footnotetext[7]{Max Planck Institute of Molecular Cell Biology and Genetics, Dresden 01307, Germany.}
\footnotetext[8]{Centre for Systems Biology Dresden, Dresden 01307, Germany.}
\footnotetext[9]{Faculty of Mathematics, Technische Universit\"{a}t Dresden, Dresden 01062, Germany.}
\footnotetext[10]{Corresponding author(s). Maria-Jose Jimenez: majiro@us.es; Heather A. Harrington: harrington@mpi.cbg.de}

\section{Background on Topological Data Analysis (TDA)}
In this section we give an overview of methods from TDA that we use in our work.
For more details on simplicial complexes and homology we refer the reader to~\cite{hatcher_algebraic_2002}.
For a survey of TDA, see~\cite{oudot_persistence_2015,edelsbrunner_persistent_2017,carlsson_persistent_2020}.

\subsection{Simplicial Complexes}
The topology of a finite point cloud $X \subset \RR^d$ is discrete, and does not necessarily reflect the geometry of the data, such as pairwise distances between points, or the presence of clusters, loops, and voids at different scales.
Instead of considering the data points themselves, we can consider the $\epsilon$-neighbourhood of the data (i.e, the union of open balls of radius $\epsilon$ centered around the data points), where $\epsilon \geq 0$ is a scale parameter.
This neighborhood, which we denote by $B_{\epsilon}(X)$, often captures the geometry of $X$ in a more meaningful way.
For concrete computations, we need a combinatorial model for $B_{\epsilon}(X)$.

Constructing combinatorial structures to model relationships between data points is not a new idea; for example, the use of proximity graphs and nearest-neighbor graphs to model pairwise relationships between data points is common in many areas of data analysis.
In our case we use \emph{simplicial complexes}, which are higher-order generalization of simple undirected graphs.
\begin{definition}
	A simplicial complex is a collection of finite sets, each of which is called a \emph{simplex} (plural \emph{simplices}), such that every non-empty subset of a simplex is also a simplex.
\end{definition}
If $\sigma$ is a simplex in a simplicial complex and $\tau \subset \sigma$, we say that $\tau$ is a \emph{face} of sigma and write $\tau \leq \sigma$.
A subset of a simplicial complex $K$ that is also a simplicial complex is a called a \emph{subcomplex} of $K$.
The \emph{dimension} of a simplex $\sigma \in K$ is given by $\dim(\sigma) = \#\sigma - 1$, that is, one less than the cardinality of $\sigma$; we say that $\sigma$ is an \emph{$n$-simplex} if $\dim(\sigma) = n$.
The dimension of $K$ is the maximum dimension of any simplex in $K$.
Thus, a simple undirected graph is a 1-dimensional simplicial complex, where the vertices are 0-simplices and the edges are 1-simplices; the vertices incident to each edge are the faces of that edge; and subgraphs correspond to subcomplexes.

A simplicial complex $K$ gives rise to a topological space $|K|$, called its \emph{geometric realization}, by identifying each $n$-simplex $\sigma \in K$ with an $n$-dimensional \emph{geometric simplex} $|\sigma|$.
An $n$-dimensional geometric simplex is the convex hull of $n+1$ affinely independent points in $\RR^d$.
To illustrate, the geometric simplices corresponding to $0$, $1$, $2$, and $3$-dimensional simplices are points, line segments, triangles, and tetrahedra respectively.
To obtain the geometric realization of $K$, these geometric simplices are glued together along their common faces as determined by the combinatorial structure of $K$.
Thus, $K$ serves as a combinatorial model for the topological space $|K|$.
Every topological space is weakly homotopy equivalent to the geometric realization of a simplicial complex, so simplicial complexes are a powerful tool for approximating the topology of spaces.

While simplicial complexes are combinatorial approximations for the topology of spaces, we also need to model the relationship between spaces.
\begin{definition}
	A function $f: L \to K$ between simplicial complexes is called a \emph{simplicial map} if $\tau \leq \sigma$ in $L$ implies that $f(\tau) \leq f(\sigma)$ in $K$.
\end{definition}
An important example is when $L$ is a subcomplex of $K$, in which case the inclusion map $L \hookrightarrow K$ is a simplicial map.
Simplicial maps are a combinatorial model for continuous maps between topological spaces, in the sense that a simplicial map $f: L \to K$ induces a continuous map $|f|: |L| \to |K|$ between the geometric realizations of $L$ and $K$.
In our work, we use simplicial maps to model how the spaces formed by points of each species glue into each other.

\subsection{\v{C}ech and Alpha Complexes}
We again consider the $\epsilon$-neighbourhood of the data, $B_{\epsilon}(X)$.
A simplicial complex that models $B_{\epsilon}(X)$ is the \emph{\Cech complex} $\Cech_{\epsilon}(X)$, which is defined by
\begin{equation}
	\Cech_{\epsilon}(X) = \left\{\sigma \subseteq X : \bigcap_{x \in \sigma} B_{\epsilon}(x) \neq \emptyset\right\}.
\end{equation}
That is, simplices in $\Cech_{\epsilon}(X)$ are subsets of the data, and a subset of points $\sigma \subseteq X$ comprises a simplex in $\Cech_{\epsilon}(X)$ if and only if the $\epsilon$-balls around the points in $\sigma$ have a common point of intersection.
The following result is a consequence of a classical result in algebraic topology known as the Nerve Lemma, and it implies that $\Cech_{\epsilon}(X)$ encodes the spatial structure of $X$ at scale $\epsilon$.
\begin{theorem}
	The geometric realization of the \Cech complex $\Cech_{\epsilon}(X)$ is homotopy equivalent to the $\epsilon$-neighbourhood $B_{\epsilon}(X)$.
\end{theorem}

Another simplicial complex that encodes the spatial structure of $X$ is the \emph{alpha complex} $\Alf_{\epsilon}(X)$~\cite{Edelsbrunner1993alpha}.
A collection of points $\sigma \subset X$ is a simplex in $\Alf_{\epsilon}(X)$ if and only if $\sigma$ has a circumsphere $S$ of radius at most $\epsilon$ such that the open ball bounded by $S$ contains no points in $X$.
The alpha complex encodes the same topological information as the \Cech complex, in the following sense.
\begin{theorem}
	For each $\epsilon \geq 0$, $\Alf_{\epsilon}(X)$ is a subcomplex of $\Cech_{\epsilon}(X)$, and the induced inclusion map of the geometric realizations $|\Alf_{\epsilon}(X)| \hookrightarrow |\Cech_{\epsilon}(X)|$ is a homotopy equivalence.
\end{theorem}
Under mild conditions on the points in $X$ and for large enough $\epsilon$, the alpha complex $\Alf_{\epsilon}(X)$ is the well-known \emph{Delaunay triangulation} $\Del(X)$ of $X$~\cite{delaunay_1934}, commonly used in computer graphics, geographical information systems, and finite-element/finite-volume analysis.
For computational purposes, the alpha complex is often preferred over the \Cech complex, as it has fewer simplices and can be computed more efficiently.

\subsection{Homology}
One way to extract topological information from a simplicial complex $K$ is through its \emph{homology}.
For each $n \geq 0$, let $C_n(K)$ denote the vector space of formal $\ZZ/2\ZZ$-linear sums of $n$-dimensional simplices in $K$, with $C_{n}(K) = 0$ if $n < 0$ or $n > \dim(K)$.
The \emph{boundary} of an $n$-dimensional simplex $\sigma \in C_n(K)$, denoted by $\partial \sigma$, is the formal sum of its $(n-1)$-dimensional faces in $C_{n-1}(K)$:
\[
	\partial \sigma \coloneq \sum_{\substack{\tau \subset \sigma\\\dim(\tau) = n-1}} \tau.
\]
By linearly extending across all simplices in $C_n(K)$, we obtain the \emph{boundary map} $\partial_n: C_n(K) \to C_{n-1}(K)$.
It is straightforward to verify that $\partial_{n} \circ \partial_{n+1} = 0$, that is, $\im(\partial_{n+1}) \subseteq \ker(\partial_{n})$.
Hence, the quotient vector space
\[
	H_n(K) \coloneq \ker(\partial_n)/ \im(\partial_{n+1})
\]
is well-defined, and is called the \emph{homology of $K$ in degree $n$}.
Informally, $H_n(K)$ is generated by $n$-dimensional ``holes" in $K$; for $n = 0, 1$ and $2$, these correspond to connected components, loops, and voids in $K$ respectively.

A crucial fact about homology is that it is \emph{functorial}: a simplicial map $f: K \to L$ induces a linear map $H_n(f): H_n(K) \to H_n(L)$, and these linear maps satisfy the following properties:
\begin{enumerate}
	\item For the identity map $\mathrm{id}_K: K \to K$, we have $H_n(\mathrm{id}_K) = \mathrm{id}_{H_n(K)}$.
	\item For any two simplicial maps $f: K \to L$ and $g: L \to M$, we have $H_n(g \circ f) = H_n(g) \circ H_n(f)$.
\end{enumerate}
This is important for our work, as the functoriality of homology allows us to track how topological features in one space map to topological features in another space, and to compare the topological features of different spaces.

\subsection{Filtrations and Persistent Homology}
Let $F_{\bullet}(X)$ denote either the \Cech or alpha filtration of $X$.
An issue with the pipeline described so far is the choice of scale parameter $\epsilon$.
This can dramatically affect the topology of $|F_{\epsilon}(X)|$, and without prior knowledge of the data, it is not clear how to choose $\epsilon$.
Therefore, instead of fixing a single value of $\epsilon$, we consider the entire family of simplicial complexes $\{F_{\epsilon}(X)\}_{\epsilon \in [0, \infty]}$, which is a multiscale representation of the geometry of $X$.
\begin{definition}[Filtration]
	Let $K$ be a simplicial complex.
	A $[0, \infty]$-indexed \emph{filtration} of $K$ is a family $\{K_t\}_{t \in [0, \infty]}$ of subcomplexes $K$ such that if $s \leq t$ then $K_s$ is a subcomplex of $K_t$, and
	\[
		K = K_{\infty} = \bigcup_{t \in [0, \infty]} K_t.
	\]
	We say that $K$ is a \emph{filtered simplicial complex} if it is equipped with a filtration.
	A \emph{subfiltration} of $K_{\bullet}$ is a family of simplicial complexes $\{L_t\}_{t \in [0, \infty]}$ such that $L_t$ is a subcomplex of $K_t$ for each $t \in [0, \infty]$.
\end{definition}
It is straightforward to verify that $\{\Cech_{\epsilon}(X)\}_{\epsilon \in [0, \infty]}$ and $\{\Alf_{\epsilon}(X)\}_{\epsilon \in [0, \infty]}$ are filtrations of the power set $\mathcal{P}(X)$ of $X$; these are referred to as the \emph{\Cech filtration} and \emph{alpha filtration} respectively.

\begin{definition}[Persistent homology]
	Let $K_{\bullet}$ be a filtration of a simplicial complex $K$.
	The \emph{persistent homology} of $K_{\bullet}$ in degree $n$ is the collection of vector spaces $\{H_n(K_t)\}_{t \in [0, \infty]}$ together with the linear maps $\{H_n(K_s) \to H_n(K_t)\}_{s \leq t}$ induced by the inclusion maps $K_s \hookrightarrow K_t$.
	We denote the persistent homology of $K_{\bullet}$ in degree $n$ by $\PH_n(K_{\bullet})$.
\end{definition}
$\PH_n(F_{\bullet}(X))$ simultaneously encodes the $n$-dimensional topological features of $X$ at each scale $\epsilon$, as well as the evolution of these topological features across scales.

\subsection{Persistence Modules and Persistence Diagrams}
Persistent homology is an example of a \emph{persistence module}, that is, a collection of vector spaces $\{V_t\}_{t \in [0, \infty]}$ equipped with linear maps $\{v_{s,t}: V_s \to V_t\}_{s \leq t}$ that satisfy the following properties.
\begin{enumerate}
	\item $v_{t,t} = \mathrm{id}_{V_t}$.
	\item $v_{r,t} = v_{s,t} \circ v_{r,s}$ whenever $r \leq s \leq t$.
\end{enumerate}
Another example of a persistence module is the \emph{interval module} $I^{[s, t)}$, defined as follows:
\[
	\left(I^{[s, t)}\right)_{r} = \begin{cases}
		\ZZ/2\ZZ & \text{if } s \leq r < t, \\
		0        & \text{otherwise},
	\end{cases}
\]
with all non-zero linear maps being the identity map on $\ZZ/2\ZZ$.
Standard results from the theory of persistence modules imply that whenever $K_{\bullet}$ is a filtration of \emph{finite} simplicial complexes, then for any $n \geq 0$ its persistent homology $\PH_n(K_{\bullet})$ decomposes uniquely as a direct sum of interval modules (see for example the main result in~\cite{crawley-boevey_decomposition_2015}).
\begin{theorem}[Interval decomposition of persistent homology]
	Let $\{K_t\}_{t \in [0, \infty]}$ be a filtration of finite simplicial complexes, and let $n \geq 0$ be fixed.
	Then there exists a unique multiset of intervals $\{[b_i, d_i)\}_{i \in I}$, each contained in $[0, \infty]$, such that the persistent homology $\PH_n(K_{\bullet})$ is isomorphic to the direct sum $\bigoplus_{i \in I} I^{[b_i, d_i)}$.
\end{theorem}
Consequently, we can visualize the persistent homology of a filtration $K_{\bullet}$ as a multiset of points $\{(b_i, d_i)\}$ in the plane.
The point $(b_i, d_i)$ corresponds to the interval $[b_i, d_i)$ in the interval decomposition of $\PH_n(K_{\bullet})$, and hence represents an $n$-dimensional topological feature that is present in $K_t$ if and only if $t \in [b_i, d_i)$.

\section{Additional Details of \texorpdfstring{$k$}{k}-chromatic Gluing Maps}
Suppose $X$ is a finite set of points in $\RR^d$ that is partitioned into subsets $X_0, \ldots, X_s$; we refer to the indices in $\{0, \ldots, s\}$ as \emph{labels}, and each $X_i$ as a \emph{species}.
We adopt the same notation as in the main manuscript; for a subset $I$ of the labels, we let $X_I$ denote the subset of points in $X$ whose labels belong to $I$:
\[
	X_I \coloneqq \bigcup_{i \in I} X_i.
\]
For each such subset $I \subseteq \{0, \ldots, s\}$, we have that $\Cech_{\bullet}(X_I)$ is a subfiltration of $\Cech_{\bullet}(X)$.
By taking the disjoint union of filtrations $\Cech_{\bullet}(X_I)$ over all subsets of labels $I$ containing $k$ different labels, we obtain the \emph{$k$-chromatic gluing map}:
\[
	q_{k, \bullet}: \bigsqcup_{\substack{I \subseteq \{0, \ldots, s\}\\\#I = k}} \Cech_{\bullet}(X_I) \to \Cech_{\bullet}(X).
\]
The $k$-chromatic gluing map yields simplicial maps at each scale, which glues together the complexes spanned by each $k$-tuple of species.
The functoriality of homology then gives rise to a family of linear maps ($\epsilon \geq 0$):
\[
	H_n(q_{k,\epsilon}): \bigoplus_{\substack{I \subseteq \{0, \ldots, s\}\\ \#I=k}} H_n(\Cech_{\epsilon}(X_I)) \to H_n(\Cech_{\epsilon}(X)).
\]
The kernel, cokernel, and image of these linear maps, over all values of $\epsilon$, assemble into persistence modules due to standard results in algebra---see for example~\cite[Lemma 9.1]{lang_algebra_2002}.
These persistence modules also have interval decompositions, and can be described by and visualized through their persistence diagrams.

\subsection{Chromatic filtrations}
As mentioned in the main text, a computational difficulty that arises when using the \Cech filtration is that the number of simplices in the filtration scales exponentially with the number of data points.
This is not a problem for the alpha filtration, but on the other hand $\Alf_{\bullet}(X_I)$ is not necessarily a subfiltration of $\Alf_{\bullet}(X)$ (see for example~\cite[Figure 2]{Natarajan2024}).
This means that the $k$-chromatic gluing maps are not well-defined for the alpha filtration.
This issue of the alpha filtration is addressed by the \emph{chromatic Delaunay triangulation}~\cite{biswas2022size}, the \emph{coupled alpha filtration}~\cite{reani_coupled_2023}, the \emph{chromatic alpha filtration}~\cite{reani_coupled_2023,Montesano2026chromatic}, and the \emph{chromatic Delaunay--\Cech filtration}~\cite{Natarajan2024}.
In our work, we use the chromatic \DelCech filtration $\DelCech_{\bullet}(X)$, which is a computationally tractable subfiltration of the \Cech filtration such that $|\DelCech_{\epsilon}(X)|$ is homotopy-equivalent to $|\Cech_{\epsilon}(X)|$ at each scale $\epsilon \in [0,\infty]$ and hence captures the same topological information as the \Cech filtration.

\subsection{Comparison with \texorpdfstring{$k$}{k}-chromatic inclusion maps}
Let $F_{\bullet}(X)$ denote any one of the \Cech, chromatic alpha, or chromatic Delaunay--\Cech filtrations.
In~\cite{Montesano2026chromatic}, the authors define the \emph{$k$-chromatic subfiltration} of $F_{\bullet}(X)$---which we denote here by $F_{\bullet}(X; k)$---as comprising those simplices in $F_{\bullet}(X)$ that contain points of at most $k$ different species.
We refer to the inclusion $i_{k, \bullet} : F_{\bullet}(X; k) \hookrightarrow F_{\bullet}(X)$ as the \emph{$k$-chromatic inclusion map}.
The $k$-chromatic gluing map factors through the $k$-chromatic inclusion map, as in the following commutative diagram:
\begin{equation}
	\begin{tikzcd}[row sep=3em]
		\bigsqcup_{\substack{I \subseteq \{0, \ldots, s\}\\\#I = k}} F_{\bullet}(X_I) \arrow{r}{q_{k, \bullet}} \arrow[two heads]{d} & F_{\bullet}(X) \\
		F_{\bullet}(X; k) \arrow[hook, "i_{k, \bullet}"']{ur} &
	\end{tikzcd}
\end{equation}
Thus both maps are closely related, and in fact coincide for $k=1$.
In general however, they capture different information about the data, and \Cref{fig: gluing vs inclusion} shows an example to illustrate these differences.
The diagrams from the $k$-chromatic gluing map are more interpretable, and capture differences that are not captured by the $k$-chromatic inclusion map, as shown in \Cref{fig: gluing vs inclusion,fig: comparison 3,fig: comparison 4}.

We remark that it is possible to consider other types of inclusions as well; for example, one could consider the inclusion $F_{\bullet}(X; k)$ into $F_{\bullet}(X; k+1)$.

\section{Additional figures for the examples from the main text}
\begin{figure}[H]
	\centering
	\includegraphics[width=0.85\linewidth]{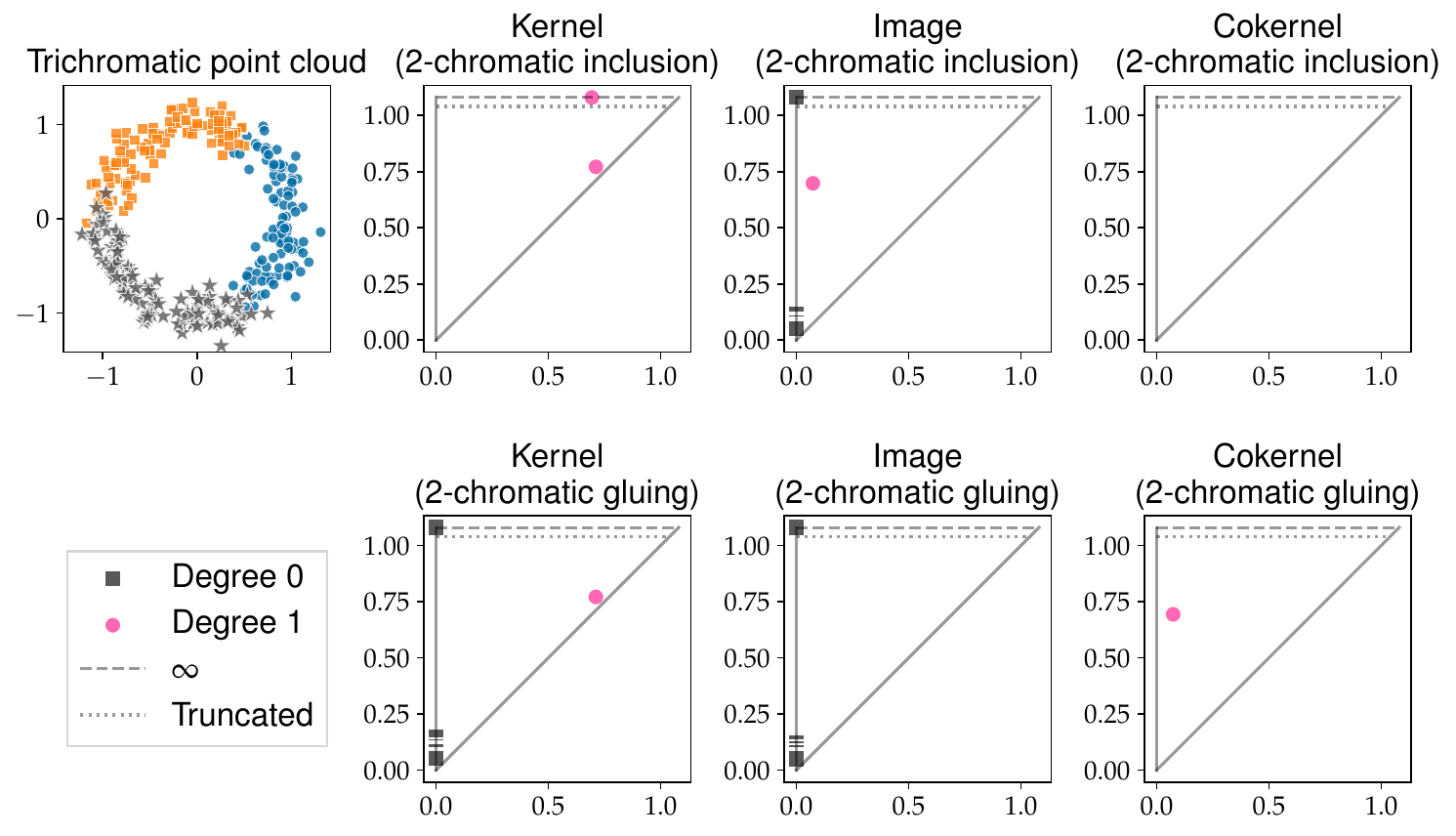}
	\caption{%
		Comparison of the kernel, image and cokernel persistence diagrams for the 2-chromatic inclusions and 2-chromatic gluing maps for a trichromatic point cloud.
		At scale $0.4$, the 2-chromatic subfiltration of the chromatic Delaunay--\Cech filtration consists of every edge in the filtration, since every edge contains points of at most 2 different species.
		Notably, it contains only a single connected component, which is the loop formed by all the points together.
		Therefore, the 2-chromatic inclusion map contains no information in the kernel in degree 0, or in the cokernel in degree 1.
		On the other hand, the disjoint union of the subfiltrations spanned by each pair of species consists of only the edges between points of the same species.
		Therefore, the 2-chromatic gluing map contains three connected components, two of which are identified with the third by the 2-chromatic gluing map.
		These two connected components give rise to two features in the kernel in degree 0 (with death time at infinity).
		The cokernel contains a cycle in degree 1, because this cycle is not present when considering any pair of species but appears when all three species are considered together.
	}
	\label{fig: gluing vs inclusion}
\end{figure}

\begin{figure}[h]
	\centering
	\def\imagewidth{0.49\linewidth}
	\def\captionwidth{0.43\linewidth}
	\begin{subfigure}[t]{\imagewidth}
		\centering
		\includegraphics[width=\linewidth]{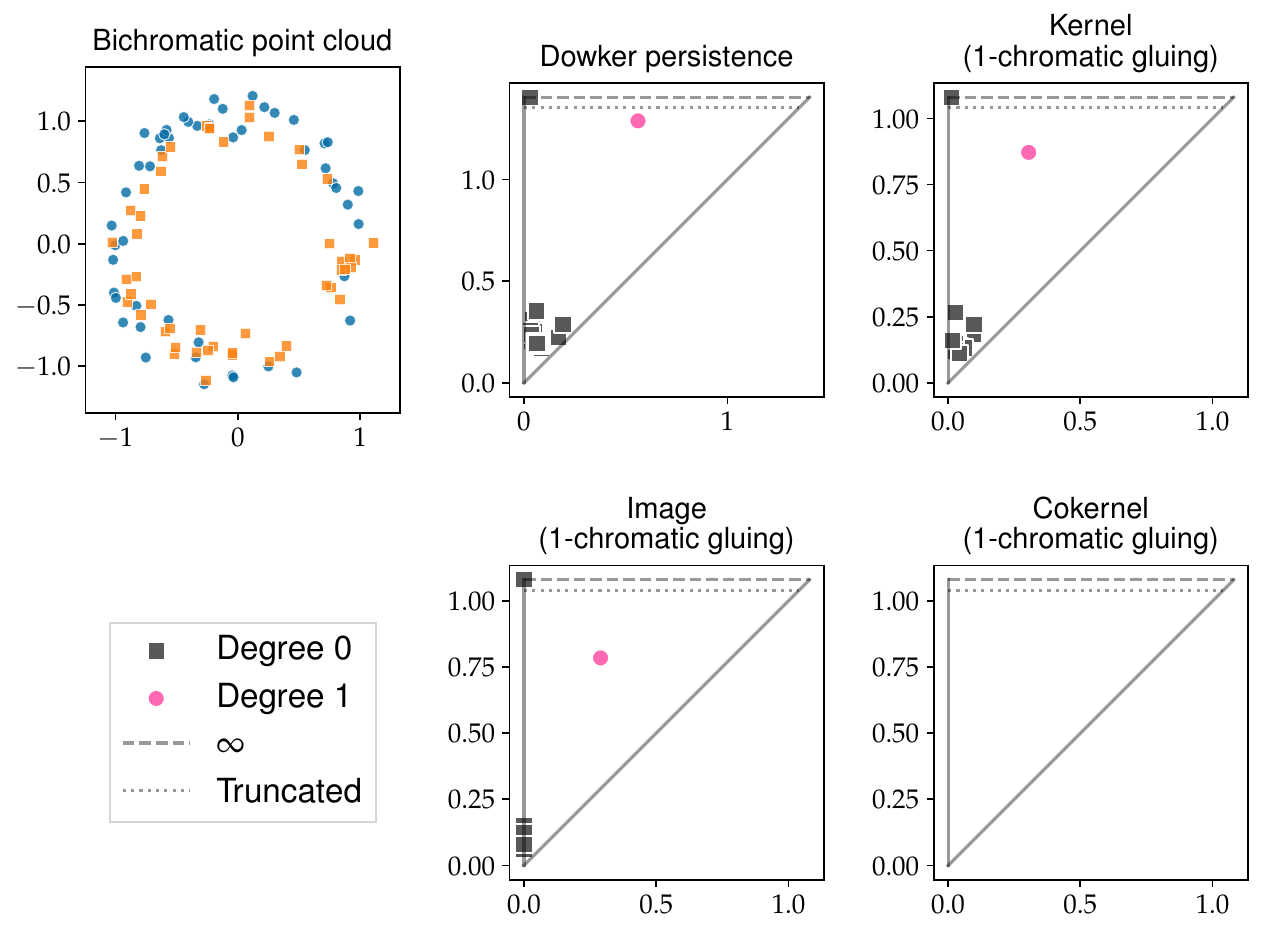}
		\caption{
		}
		\label{fig: comparison_1_1}
	\end{subfigure}\hfill
	\begin{subfigure}[t]{\imagewidth}
		\centering
		\includegraphics[ width=\linewidth]{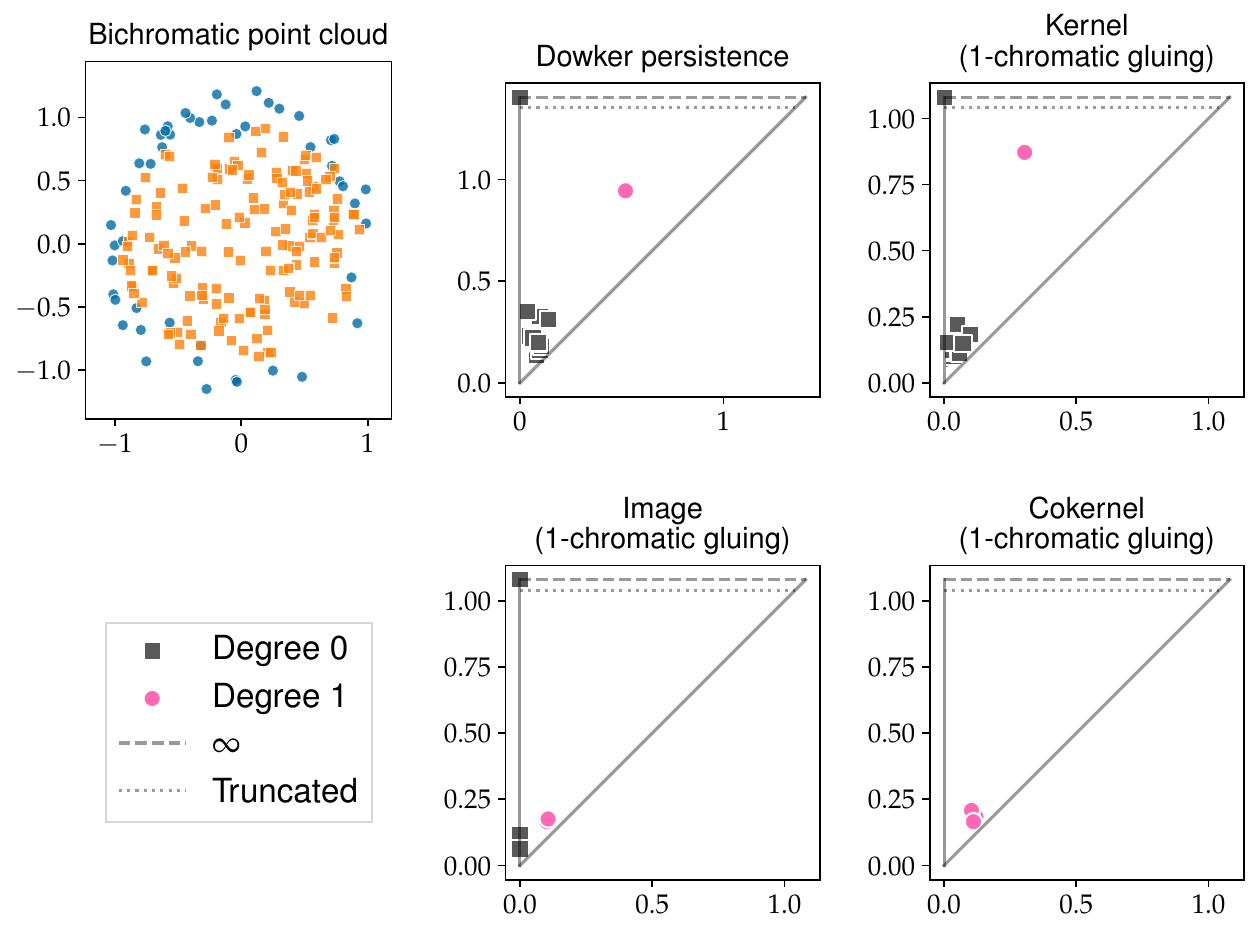}
		\caption{
		}
		\label{fig: comparison_1_2}
	\end{subfigure}\\[\baselineskip]
	\begin{subfigure}[t]{\imagewidth}
		\centering
		\includegraphics[width=\linewidth]{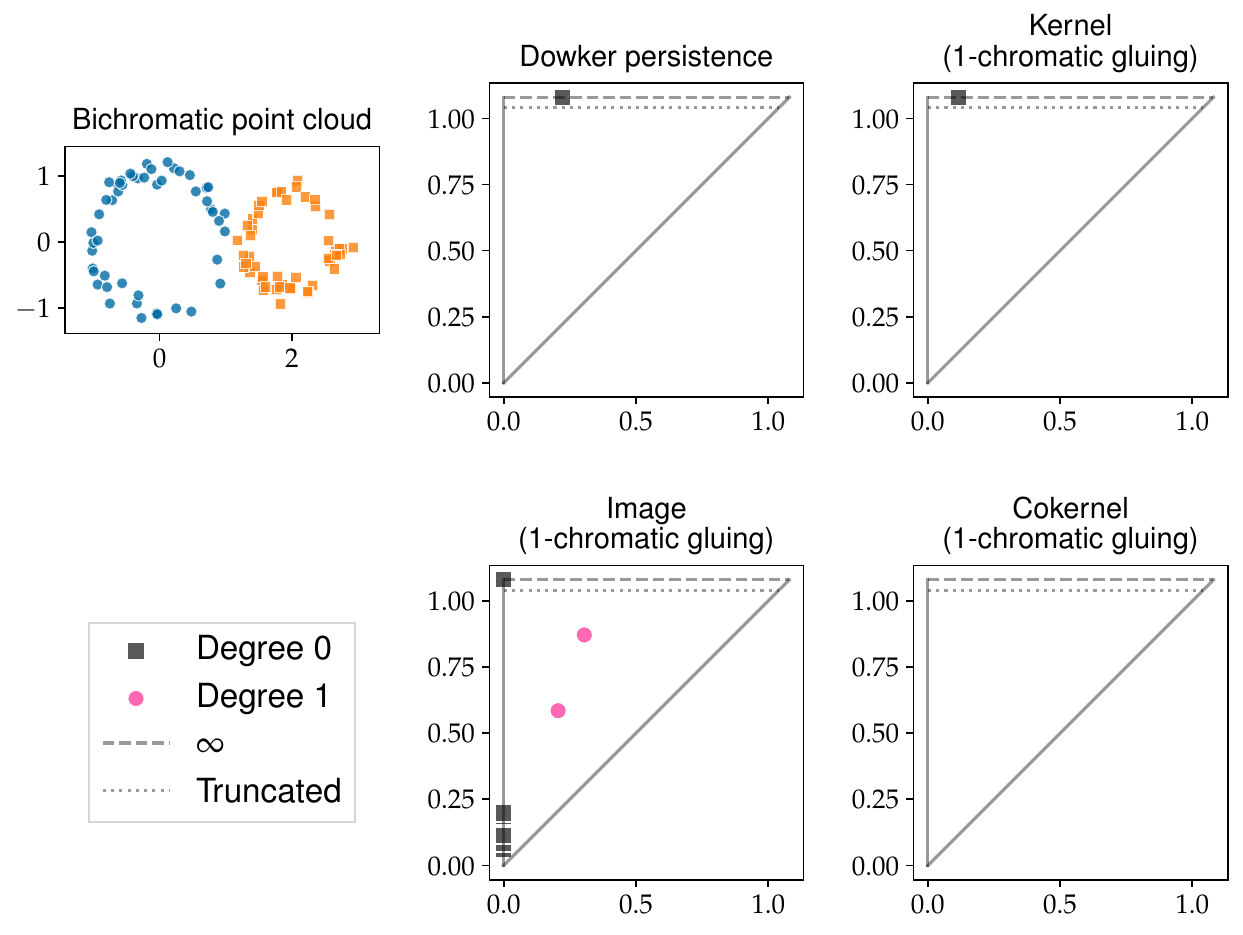}
		\caption{
		}
		\label{fig: comparison_2_1}
	\end{subfigure}\hfill
	\begin{subfigure}[t]{\imagewidth}
		\centering
		\includegraphics[width=\linewidth]{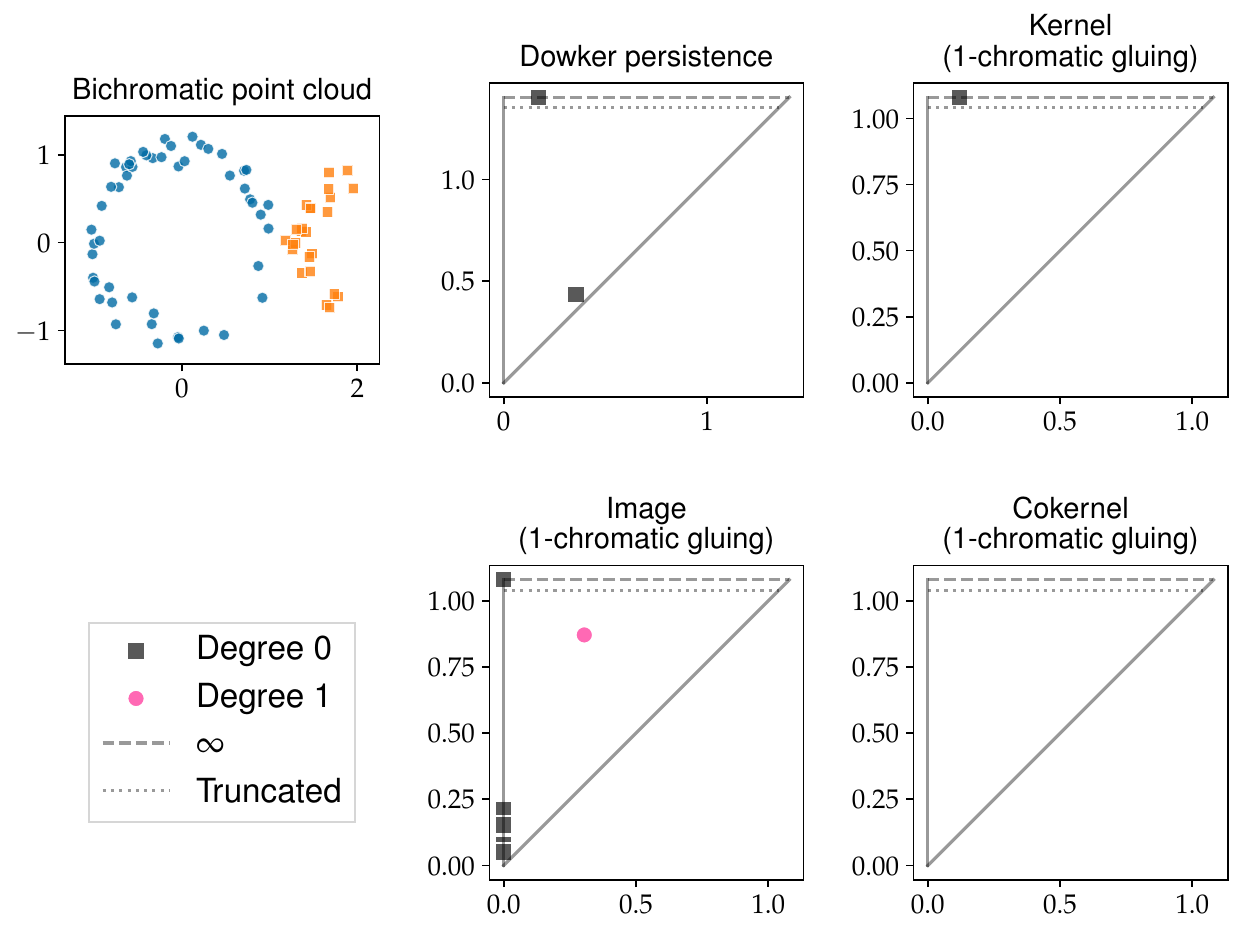}
		\caption{
		}
		\label{fig: comparison_2_2}
	\end{subfigure}
	\caption{%
		Diagrams corresponding to Dowker persistence, as well as kernel, image and cokernel for 1-chromatic gluing to illustrate the two first rows in the table in Figure 1(c) in the main manuscript.
		Examples (a) and (b) (top left and top right) can be distinguished by the persistence diagrams from the 1-chromatic gluing map, but not by the corresponding Dowker persistence diagrams.
		Similarly, examples (c) and (d) (bottom left and bottom right) can be distinguished by the 1-chromatic gluing map, but not by the corresponding Dowker persistence diagrams.
	}
	\label{fig: comparison 1 and 2}
\end{figure}

\begin{enumerate}
	\item \Cref{fig: gluing vs inclusion} shows an example to highlight the difference between $k$-chromatic gluing map and $k$-chromatic inclusion map.
	\item \Cref{fig: comparison 1 and 2} compares our framework with Dowker persistence, expanding on the examples in the first two rows of the table in Figure 1(c) in the main text.
	\item \Cref{fig: comparison 3,fig: comparison 4} expand on the examples in the third and fourth rows of the table in Figure 1(c) in the main text, comparing the persistence diagrams obtained from the $k$-chromatic gluing map with those from the $k$-chromatic subfiltrations.
	\item \Cref{fig: examples-SI} expands upon the examples from Figure 2 in the main text, showing the remaining persistence diagrams for some of those examples.
\end{enumerate}

\begin{figure}[h]
	\centering
	\begin{subfigure}{\linewidth}
		\centering
		\includegraphics[width=0.9\linewidth]{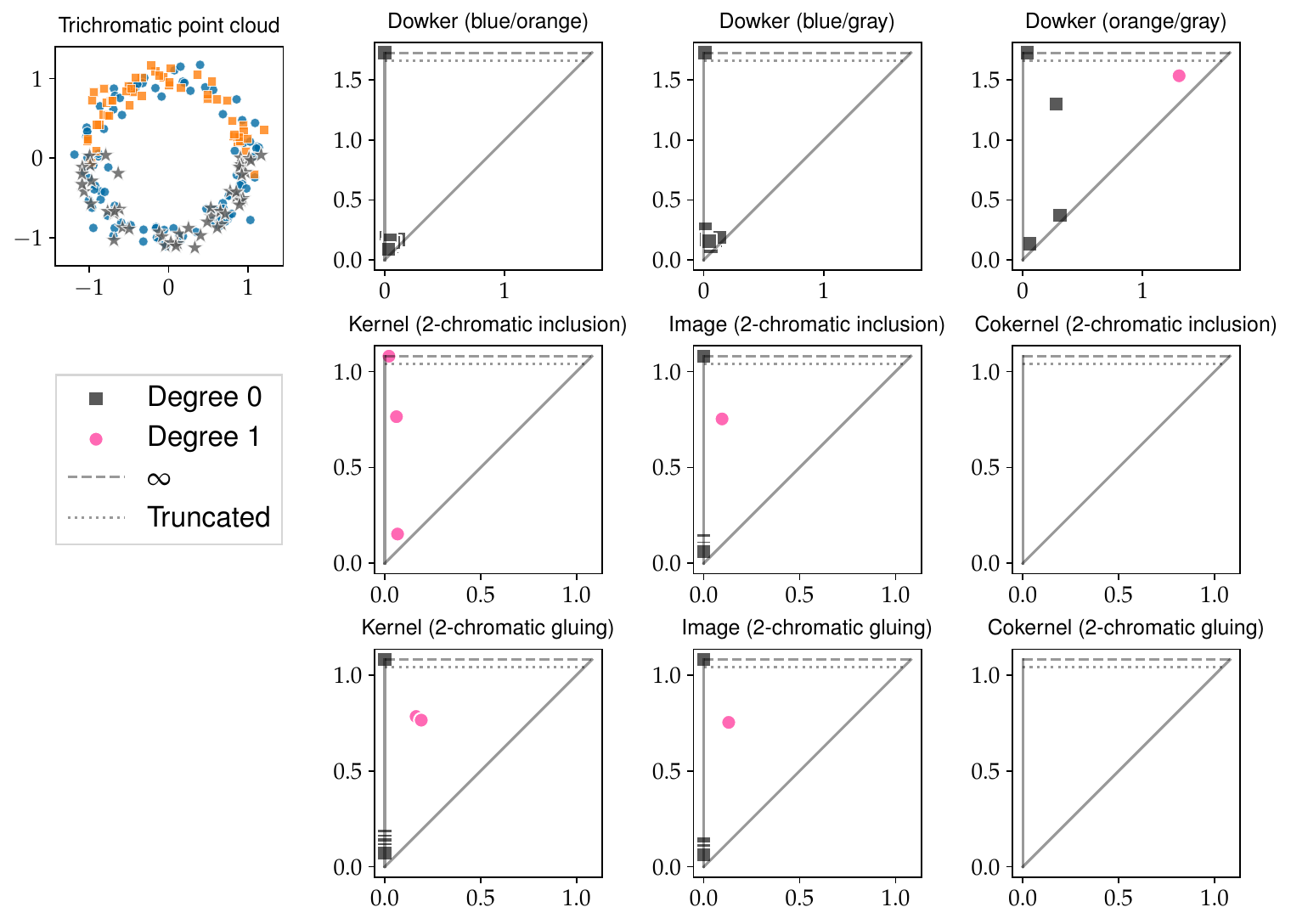}
		\vspace{-0.15cm}
		\caption{Each of the three combinations of two colours forms a loop, which is not clearly seen in any of the Dowker diagrams. Two of those loops are identified with the third by the 2-chromatic gluing map when all three colors are considered together.
			This identification gives rise to the two prominent features in the kernel in degree 1 (bottom row left), and the presence of the cycle in the image diagram (bottom row middle).
			While the image diagram for the 2-chromatic inclusion map (middle row, middle) also contains a prominent feature in degree 1, the kernel contains a feature with infinite lifetime in degree 1.
			This feature corresponds to a loop formed by a path of edges of at most 2 colors in the filtration that is never filled in by triangles of at most 2 colors.
			This makes its interpretation less clear.
			It is also not clear how to recover the fact that there are multiple coincident loops in the data, each formed by 2 colors, from the diagrams of the 2-chromatic inclusion map.
		}
		\label{fig: comparison3_1}
	\end{subfigure}
\end{figure}
\begin{figure}[h]
	\ContinuedFloat
	\centering
	\begin{subfigure}{\linewidth}
		\centering
		\includegraphics[width=0.9\linewidth]{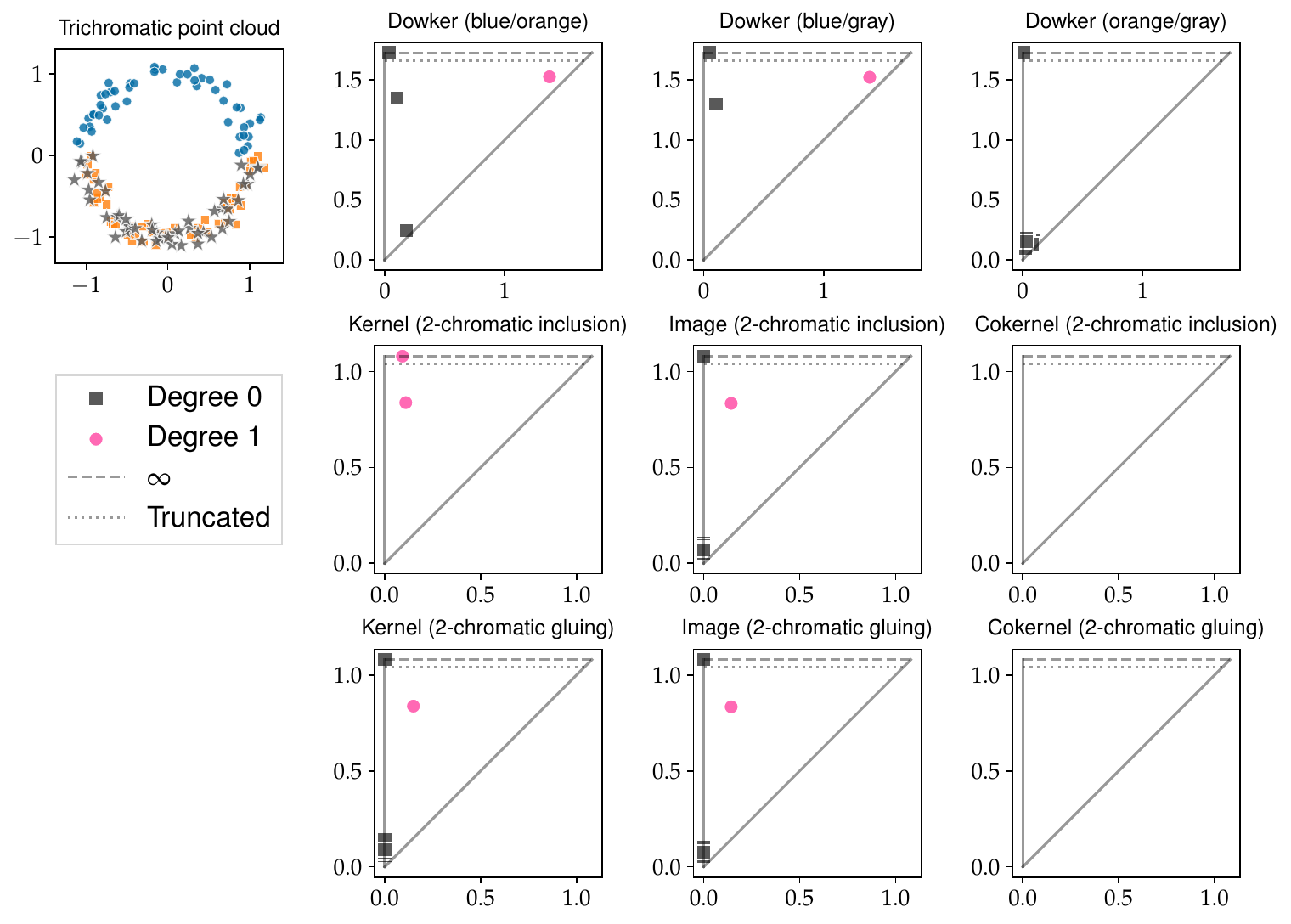}
		\vspace{-0.15cm}
		\caption{%
			In this example, the orange/blue and gray/blue (square/circle and star/circle respectively) pairs of species each form a loop, which can not be seen clearly in the Dowker diagrams.
			When all three colours are considered together, these two loops are identified 
			by the 2-chromatic gluing map.
			This results in a single prominent feature in the kernel in degree 1 (bottom row left), and a single prominent feature in the image diagram (bottom row middle).
			The interpretation of the diagrams from the 2-chromatic inclusion map (top row) is similar to the previous example (a).
		}
		\label{fig: comparison3_2}
	\end{subfigure}\\[\baselineskip]
	\vspace{-0.5cm}
	\caption{%
		These examples illustrate the third row in the table in Figure 1(c) in the main manuscript.
		The two point cloud configurations shown here can be distinguished by either the three Dowker diagrams of the three pair combinations or by the kernel diagram of 2-chromatic gluing map, but not by the diagrams from the 2-chromatic inclusion map.
		It is also easier to interpret the information from the diagrams of the 2-chromatic gluing map, as each of the features in these diagrams corresponds directly to a topological feature in some combination of colors.
	}
	\label{fig: comparison 3}
\end{figure}

\begin{figure}[h]
	\centering
	\begin{subfigure}{\linewidth}
		\centering
		\includegraphics[width=0.9\linewidth]{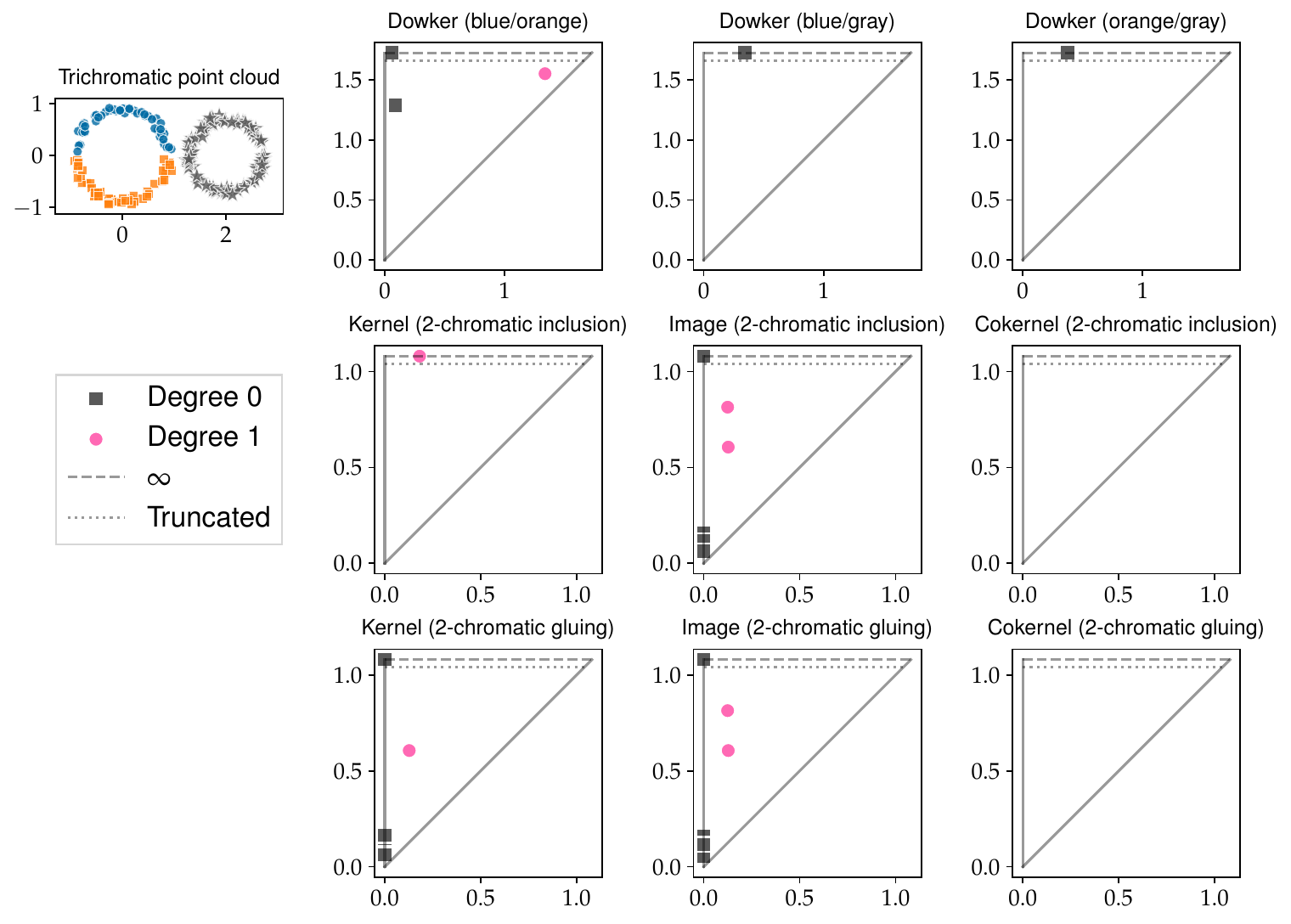}
		\vspace{-0.15cm}
		\caption{%
			There are two prominent loops in the data, one of which is formed by the orange/blue (square/circle) pair of species, while the other loop is formed by the gray (star) species. None of the Dowker diagrams captures well this information.
			The gray loop can also be seen in the gray/blue and gray/orange combinations, and these are identified with each other by the 2-chromatic gluing map when all three colors are considered together.
			This results in a prominent feature in the kernel in degree 1 (bottom row left).
			The gray loop and the orange/blue loops are both present as prominent features in the image diagram from both the 2-chromatic gluing (bottom row middle) and the 2-chromatic inclusion (top row middle) maps.
		}
		\label{fig: comparison4_1}
	\end{subfigure}
\end{figure}
\begin{figure}[h]
	\ContinuedFloat
	\centering
	\begin{subfigure}{\linewidth}
		\centering
		\includegraphics[width=0.9\linewidth]{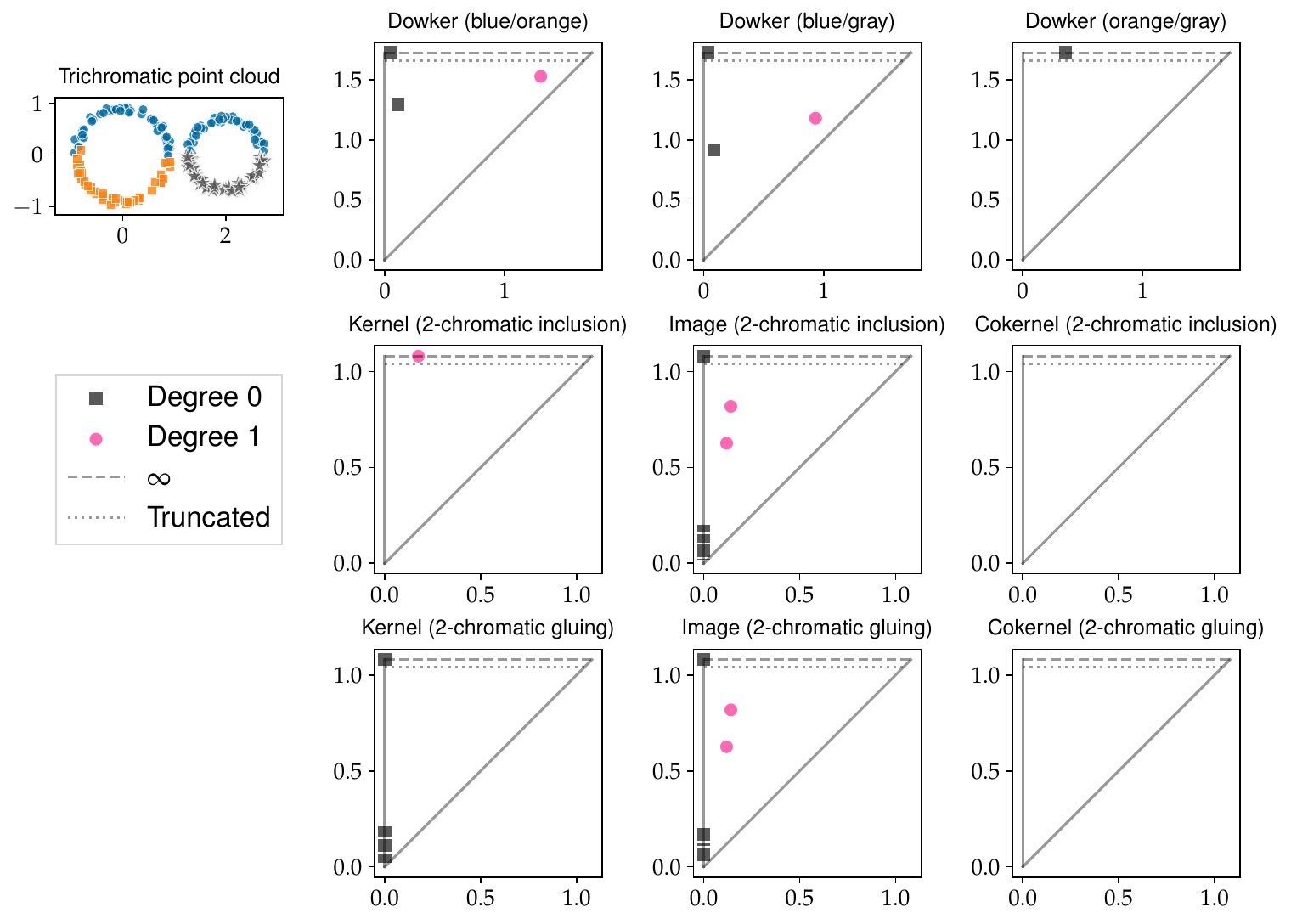}
		\vspace{-0.15cm}
		\caption{%
			This example differs from (a) in that the gray loop is replaced by a gray/blue loop, that is, there is no longer a prominent loop in the gray/orange pair of colors.
			Therefore, there is no longer a prominent feature in the kernel of the 2-chromatic gluing map in degree 1 (bottom row left).
		}
		\label{fig: comparison4_2}
	\end{subfigure}\\[\baselineskip]
	\vspace{-0.5cm}
	\caption{%
		These examples illustrate the fourth row in the table in Figure 1(c) in the main manuscript.
		The two point cloud configurations shown here can be distinguished either by the three Dowker diagrams of the three pairs of colors or by the kernel of the 2-chromatic gluing map, but not by the diagrams from the 2-chromatic inclusion map.
		As in \Cref{fig: comparison 3}, it is easier to interpret the information from the diagrams of the 2-chromatic gluing map, as each of the features in these diagrams corresponds directly to a topological feature in some combination of colors.
	}
	\label{fig: comparison 4}
\end{figure}

\begin{figure}[h]
	\centering
	\def\imagewidth{0.49\linewidth}
	\def\captionwidth{0.43\linewidth}
	\begin{subfigure}[t]{\imagewidth}
		\centering
		\includegraphics[ width=\linewidth]{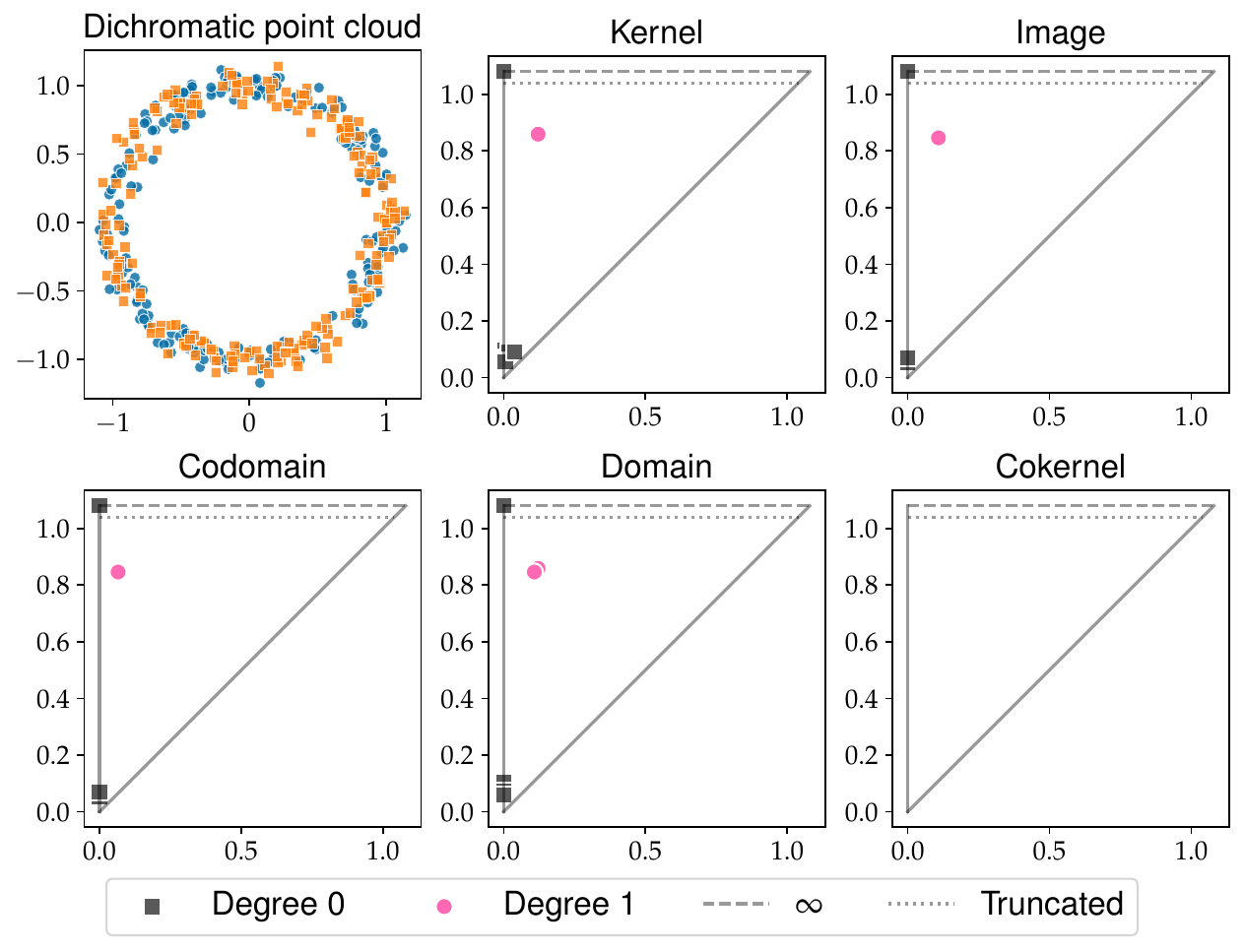}
		\caption{
		}
		\label{fig: example2b-SI}
	\end{subfigure}\hfill
	\begin{subfigure}[t]{\imagewidth}
		\centering
		\includegraphics[ width=\linewidth]{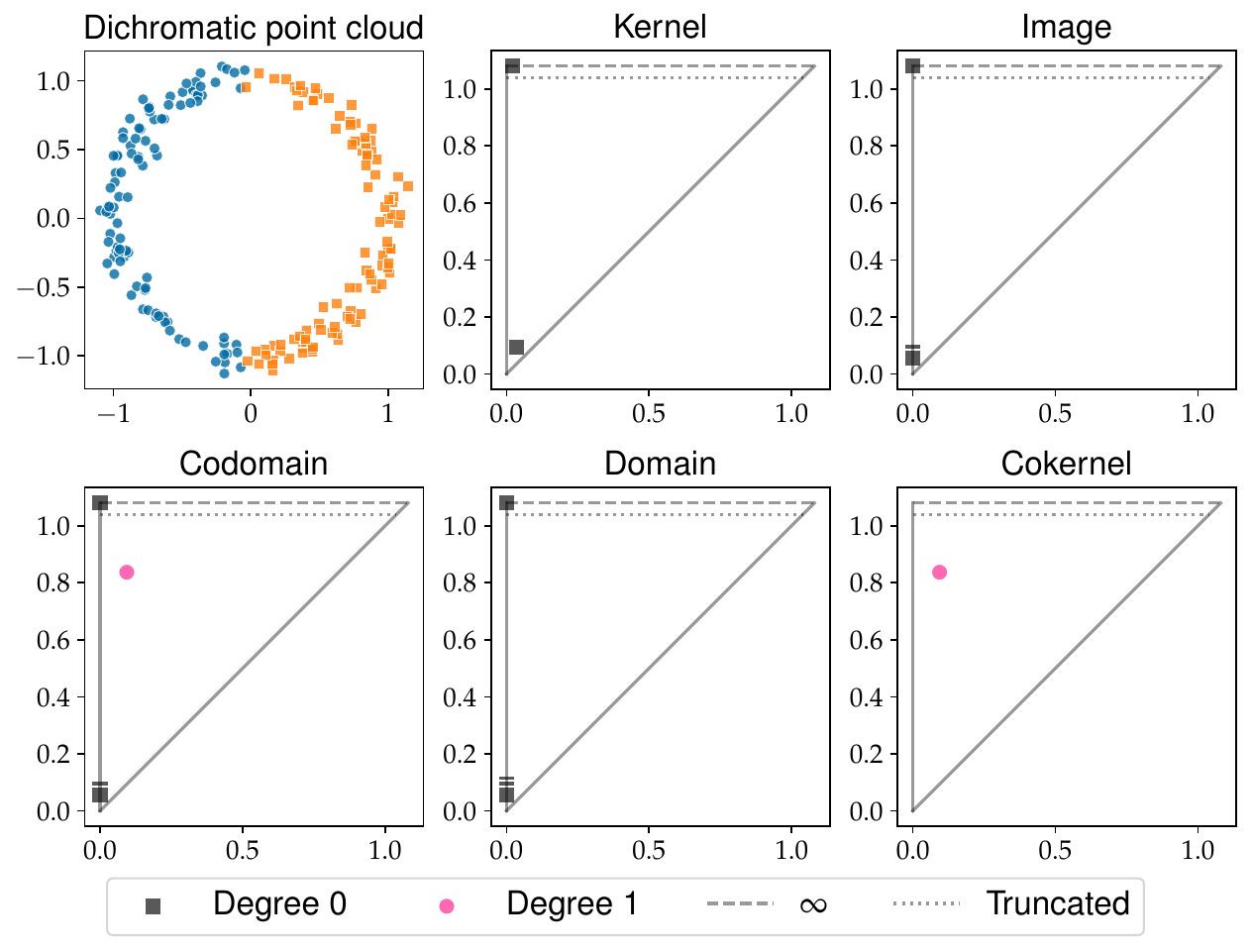}
		\caption{
		}
		\label{fig: example2c-SI}
	\end{subfigure}\\[\baselineskip]
	\begin{subfigure}[t]{\imagewidth}
		\centering
		\includegraphics[width=\linewidth]{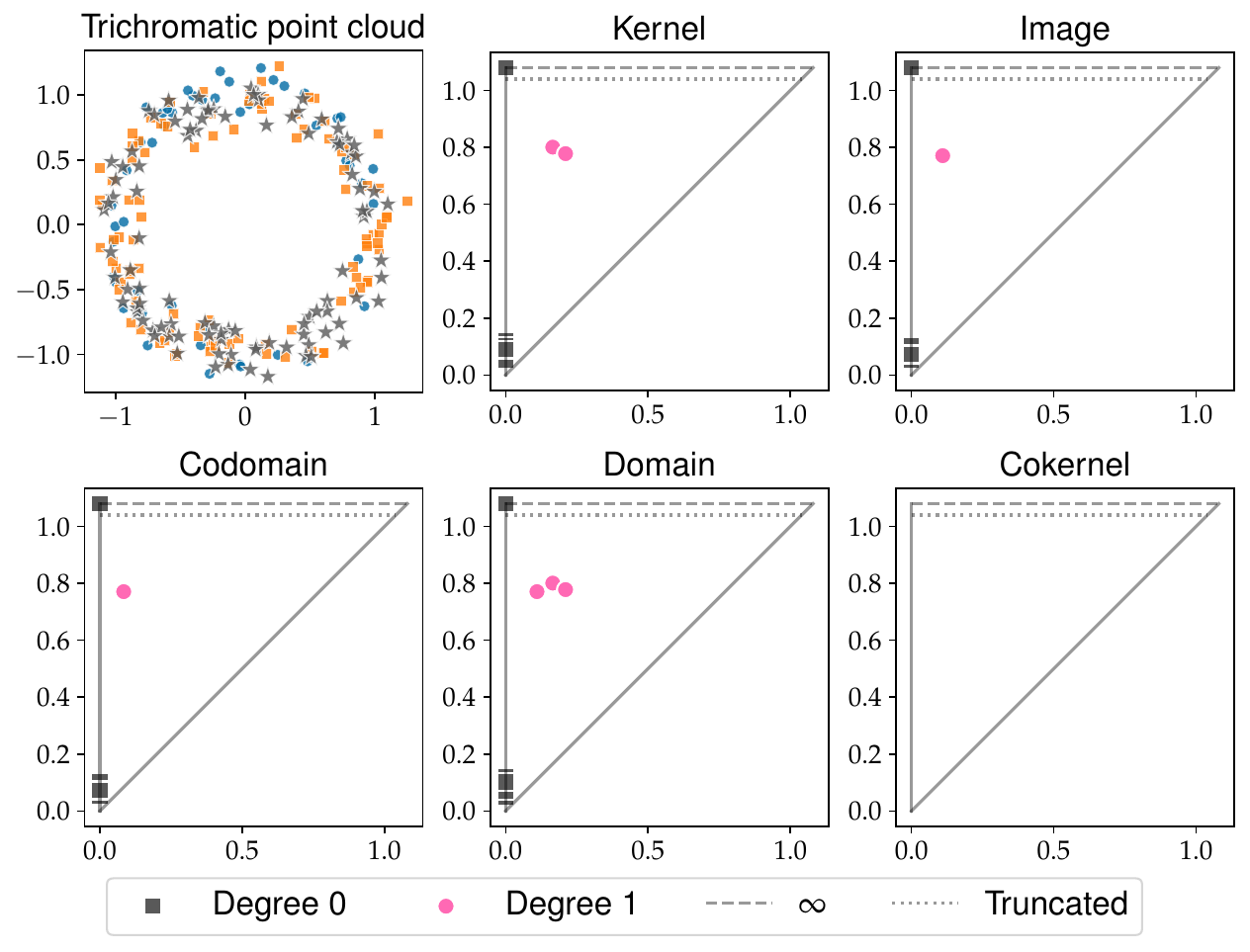}
		\caption{
		}
		\label{fig: example3-SI}
	\end{subfigure}\hfill
	\begin{subfigure}[t]{\imagewidth}
		\centering
		\includegraphics[width=\linewidth]{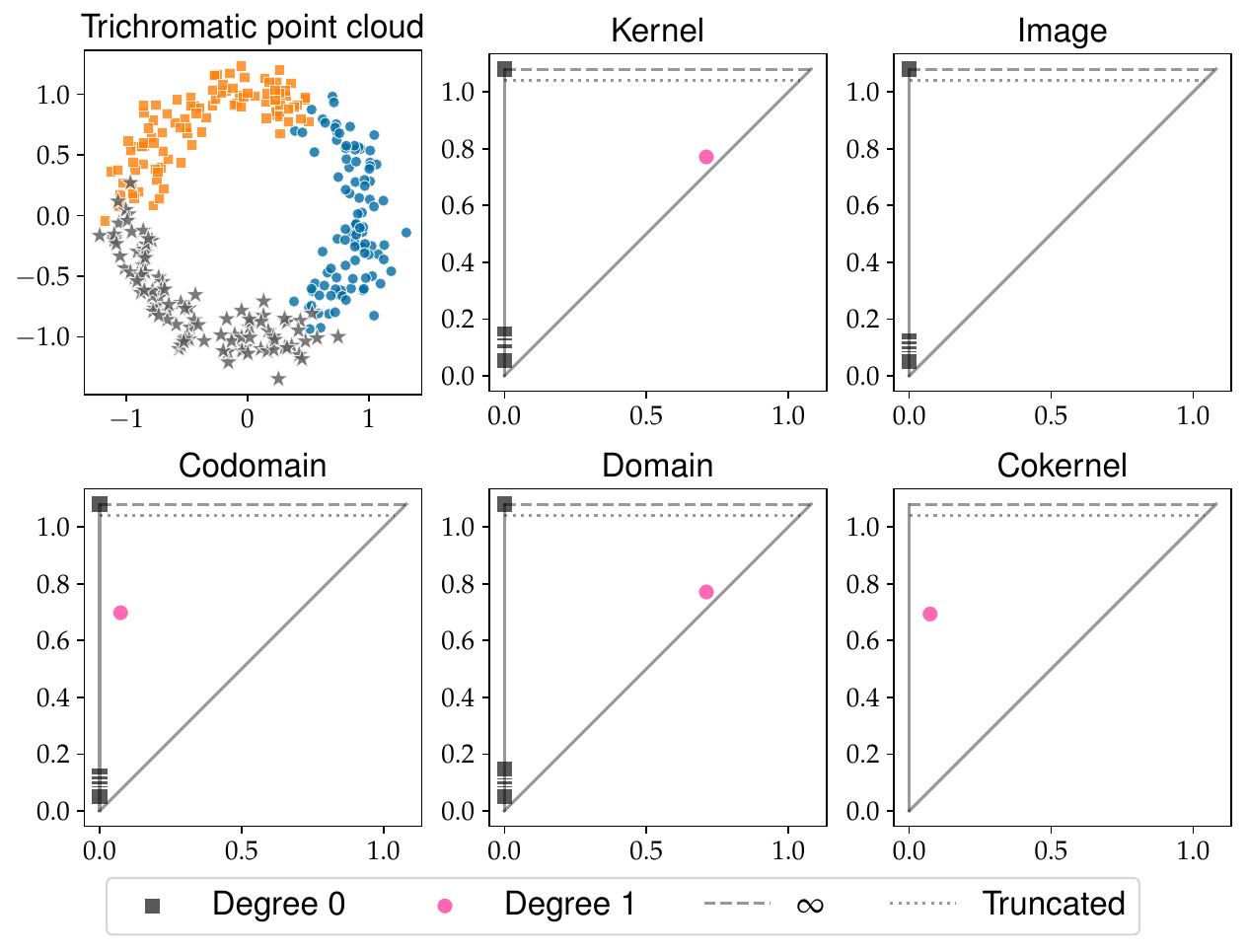}
		\caption{
		}
		\label{fig: example4-SI}
	\end{subfigure}
	\caption{%
		Additional persistence diagrams for the examples from Figure 2 of the main text.
		Examples (a) and (b) (top left and right respectively) show dichromatic point clouds, and the diagrams from the associated 1-chromatic gluing maps.
		Examples (c) and (d) (bottom left and right respectively) show trichromatic point clouds and the diagrams from the associated 2-chromatic gluing maps.
		When multiple loops are identified with each other in the full set of colors, one of them is visible in the image diagram, while the rest lie in the kernel diagram (see the 2 colors case in (a) and the three colors case in (c)).
		The cokernel diagram of the 1-chromatic gluing map captures loops formed by two colors but not by one color in (b), and the cokernel of the 2-chromatic gluing map captures the loop formed by 3 colors but not by 2 colors in (d).
	} \label{fig: examples-SI}
\end{figure}

\clearpage

\printbibliography